\documentclass[12pt,letterpaper]{article}
\usepackage{amsmath,amssymb,amsfonts,amsthm,mathtools}
\usepackage{graphicx}
\usepackage{enumerate}
\usepackage[authoryear,round]{natbib}
\usepackage{subcaption}
\usepackage[margin=1in]{geometry}
\usepackage[colorlinks=true,allcolors=black,bookmarks=false,menucolor=black]{hyperref}
\usepackage[title]{appendix}

\newcommand{\calZ}{\mathcal{Z}}
\newcommand{\bfa}{\mathbf{a}}
\newcommand{\bfz}{\mathbf{z}}
\newcommand{\bfZ}{\mathbf{Z}}
\newcommand{\bbN}{\mathbb{N}}
\newcommand{\bbR}{\mathbb{R}}
\newcommand{\bbE}{\mathbb{E}}
\newcommand{\calP}{\mathcal{P}}
\newcommand{\calB}{\mathcal{B}}
\newcommand{\rmd}{\mathrm{d}}
\newcommand{\eps}{\varepsilon}

\newcommand{\KL}{\operatorname{KL}}
\newcommand{\TTT}{\mathcal{T}}
\newcommand{\bfx}{\mathbf{x}}
\newcommand{\bfQ}{\mathbf{Q}}
\newcommand{\bbP}{\mathbb{P}}

\newcommand{\figref}[1]{Figure~\ref{#1}}
\DeclareMathOperator*{\given}{{\mid}}

\newtheorem{thm}{Theorem}
\newtheorem{lem}{Lemma}
\newtheorem{cor}{Corollary}
\theoremstyle{definition}
\newtheorem{rmk}{Remark}

\allowdisplaybreaks

\begin{document}


\title{
	Asymptotic guarantees for Bayesian phylogenetic tree reconstruction
}
\author{
	Alisa Kirichenko\thanks{Corresponding author: \texttt{alicekirichenko@optiver.com}.}\\
	Department of Statistics, University of Warwick
	\and
	Luke J. Kelly\thanks{Supported by the French government under management of Agence Nationale de la Recherche as part of the ABSint (ANR-18-CE40-0034) and PRAIRIE (ANR-19-P3IA-0001) programs.} \\
	School of Mathematical Sciences, University College Cork
	\and
	Jere Koskela\thanks{Funded by UKRI through EPSRC grant EP/V049208/1.}\\
	School of Mathematics, Statistics and Physics, Newcastle University, \\ and Department of Statistics, University of Warwick
}
\date{}
\maketitle

\begin{abstract}
We derive tractable criteria for the consistency of Bayesian tree reconstruction procedures, which constitute a central class of algorithms for inferring common ancestry among DNA sequence samples in phylogenetics.
Our results encompass several Bayesian algorithms in widespread use, such as \texttt{BEAST}, \texttt{MrBayes}, and \texttt{RevBayes}.
Unlike essentially all existing asymptotic guarantees for tree reconstruction, we require no discretization or boundedness assumptions on branch lengths.
Our results are also very flexible, and easy to adapt to variations of the underlying inference problem. We demonstrate the practicality of our criteria on two examples: a Kingman coalescent prior on rooted, ultrametric trees, and an independence prior on unconstrained binary trees, though we emphasize that our result also applies to non-binary tree models.
In both cases, the convergence rate we obtain matches known, frequentist results obtained using stronger boundedness assumptions, up to logarithmic factors.
\end{abstract}

\noindent%

\section{Introduction}
\label{sec:intro}
Reconstructing the genealogical tree of a sample of DNA sequences is a central task in phylogenetics.
Genealogies cannot usually be observed directly, but they are informative about the evolutionary history of the population and encode the positive correlations among the samples brought about by common ancestry.
There are many well-established methods for constructing a genealogy from a sample of sequences.
Some recent contributions include parsimony- or consensus-based methods \citep{sayyari/mirarab:2016, zhangetal:2018}, maximum likelihood methods which seek to maximize the likelihood of the observed sample under a given generative model of DNA sequence diversity, possibly subject to regularization \citep{stamatakis:2006, zhangetal:2021}, and Bayesian methods which seek to compute (or, more practically, sample from) the posterior distribution of genealogical trees given sequences observed at their leaves \citep{borges/kosiol:2020, ronquistetal:2012, suchardetal:2018}.
For introductions to mathematical and computational phylogenetics, see \cite{semple/steel:2013, warnow:2017}.

Our focus is consistency of Bayesian tree reconstruction procedures.
Informally, a procedure is consistent if, in the absence of model error, it recovers the true, data-generating genealogical tree with certainty in an idealized, infinite-data limit, which we take to be the length of the observed sequences.
Consistency is widely regarded as a minimal reliability criterion for Bayesian inference, and general conditions for it have been established: see, for example, \cite{schwartz:1965}.
When consistency fails, pathological results such as posterior distributions which express ever-increasing confidence in an incorrect model can arise \citep{freedman/diaconis:1983}.

Consistency of tree reconstruction has been established for many consensus methods \cite[and references therein]{markiv/eulenstein:2021}, as well as maximum likelihood methods with \citep{dinhetal:2018, zhangetal:2021} and without \citep{roch/sly:2017} regularization.
These results rest on a variety of strong assumptions: \cite{roch/sly:2017} assume that branch lengths take values in a discrete, finite set which is known a priori, \cite{dinhetal:2018} require branch lengths to be bounded above and away from zero, while \cite{zhangetal:2021} assume branch lengths are only bounded from above.
Consistency in the Bayesian setting has also been established, but only for the tree topology, treating branch lengths as nuisance parameters \citep{steel:2013}.

Our contribution is to establish Bayesian posterior consistency jointly for the tree topology and its branch lengths, regarding branch lengths as quantities of interest rather than nuisance parameters.
Our approach requires no boundedness or discretization assumptions.
Moreover, we demonstrate that phylogenetic tree models are amenable to standard tools of Bayesian asymptotic theory \citep{ghosaletal:2000}, whereas the results of \cite{steel:2013} were established using bespoke arguments for genealogical trees.
Our method is highly flexible, and can be easily applied to, for example, tree models with non-binary branching and hence a variable number of edges, and to constrained models such as ultrametric trees.
Our results apply directly to several prominent Bayesian tree inference packages: \texttt{MrBayes} \citep{ronquistetal:2012}, \texttt{BEAST 2} \citep{bouckaertetal:2019}, and \texttt{RevBayes} \citep{hohna2016revbayes}.

We demonstrate the utility of our result on two examples: unrooted, binary trees with independent branch lengths, and rooted, ultrametric binary trees with a Kingman coalescent prior \citep{Kingman82a}.
In both cases we obtain a convergence rate of $( \log k )^{ 3/4 + \beta / 2 } / \sqrt{k}$ for any $\beta > 1$, where $k$ is the length of the observed DNA sequences, matching known rates for maximum likelihood procedures up to logarithmic corrections \cite[Theorem 5.2]{dinhetal:2018}.

The rest of the article is structured as follows.
In Section \ref{section:preliminaries} we recall the general Bayesian posterior consistency theorem of \citet[Theorem 2.1]{ghosaletal:2000} which will serve as our main mathematical tool.
We will also introduce the state space and likelihood model for our genealogical trees.
In Section \ref{section:results} we state and prove our main results by providing tractable consistency criteria on Bayesian tree priors.
Section \ref{section:examples} concludes by demonstrating that two examples of standard priors satisfy our criteria, and Section \ref{section:discussion} concludes with a discussion.
Proofs of results are postponed to Sections \ref{section:proofs} and \ref{section:technical_proofs}.

\section{Preliminaries}\label{section:preliminaries}

\subsection{Bayesian posterior consistency}

Let $\mathcal{P}=\{ P_{ \theta } \}_{ \theta \in \Theta }$ be a collection of probability distributions with common support $\calZ$ and corresponding probability densities $\{ p_{ \theta } \}_{ \theta \in \Theta }$ relative to a given dominating measure.
Let $\bfZ^{(k)}=(\bfZ_1,\dots,\bfZ_k)$ be a vector consisting of independent and identically distributed random variables with common distribution $P_{\theta_0}$ for $\theta_0 \in \Theta$.
Finally, let $\bfz^{(k)}=(\bfz_1,\dots,\bfz_k)$ denote a vector arising as the observed value of $\bfZ^{(k)}$, and let $\Pi$ be a prior distribution on $\calP$.
Under appropriate identifiability conditions, we can always define a prior $\Pi$ on $\calP$ using a prior $\tilde\Pi$ on the set $\Theta$, and push it forward to $\calP$ via the mapping $\theta \mapsto P_{ \theta }$.
Throughout the paper we will slightly abuse the notation by using $\Pi$ to denote both priors.

The central object in Bayesian inference is the posterior distribution, defined for $A \in \calB( \Theta )$ as
\begin{equation*}
\Pi( A | \bfz^{(k)} ) := \frac{ \int_A p_{ \theta }( \bfz^{(k)} ) \Pi( \rmd \theta ) }{ \int_{ \Theta } p_{ \theta }( \bfz^{(k)} ) \Pi( \rmd \theta ) }.
\end{equation*}
For a sequence $\eps_k \to 0$, the posterior is said to contract around the truth $P_{\theta_0}$ at a rate $\eps_k$ if
\begin{equation}
\Pi( d_H( P_{ \theta }, P_{ \theta_0 } ) > M_k\eps_k \given \bfZ^{(k)} ) \to 0
\end{equation}
as $k \to \infty$ in $P_{ \theta_0 }$-probability for every sequence $M_k \to \infty$, where $d_H$ is the Hellinger distance,
\begin{equation*}
d_H( P_{ \theta }, P_{ \theta_0 } ) := \frac{ 1 }{ 2 } \int_{ \calZ } \Big( \sqrt{ p_{ \theta }( \bfz ) } - \sqrt{ p_{ \theta_0 }( \bfz ) } \Big)^2 \rmd \bfz.
\end{equation*}

In order to establish posterior consistency criteria and derive contraction rates for phylogenetic tree priors, we are going to employ the general machinery first introduced by \cite{ghosaletal:2000}.
More specifically, we use Theorem 8.9 from \cite{gvdv}. Consider the following statistical distances:
for $P_{\theta_0}, P_\theta \in \calP$, let
\begin{align*}
\KL( P_{\theta_0}, P_\theta ) &:= \int_{ \calZ } \log\Bigg( \frac{ p_{ \theta_0 }( \bfz ) }{ p_{ \theta }( \bfz ) } \Bigg) P_{ \theta_0 }( \rmd \bfz ), \\
V_0( P_{\theta_0}, P_\theta ) &:= \int_{ \calZ } \log^2\Bigg( \frac{ p_{ \theta_0 }( \bfz ) }{ p_{ \theta }( \bfz ) } \Bigg) P_{ \theta_0 }( \rmd \bfz ) - [ \KL( P_{\theta_0}, P_{\theta} ) ]^2
\end{align*}
be, respectively, the Kullback--Leibler divergence and the centered second Kullback--Leibler variation.
Additionally, for a set $A\subseteq \calP$ and $\eps>0$ let $N(\eps,A,d_H)$ be the $\eps$-covering number of $A$ with respect to $d_H$.
Then the following result holds.
\begin{thm}[Theorem 8.9 of \cite{gvdv}]
\label{thm:ghosaletal:2000}
Suppose there exist sets $\calP_k\subseteq \calP$, sequences $\eps_k \to 0$ and $\zeta_k\to 0$ such that $k\zeta_k^2\to \infty$ and $\zeta_k\leq\eps_k$, as well as a constant $C>0$ such that
\begin{align}
\label{eq:remaining}
\Pi(\calP \backslash \calP_k) &\leq \exp( - ( C + 4 ) k \zeta_k^2 ), \\
\label{eq:entropy}
\log N(\eps_k/2,\calP_k, d_H)&\leq k\eps_k^2, \\
\label{eq:prior}
\Pi( P_\theta : \KL( P_{\theta_0}, P_\theta )\leq \zeta_k^2, V_0( P_{\theta_0},P_\theta )\leq \zeta_k^2 ) &\geq \exp( - C k \zeta_k^2 ).
\end{align}
Then the posterior $\Pi( \cdot | \bfz^{(k)} )$ converges to the data-generating distribution $P_{ \theta_0 }$ at rate $\eps_k$.
\end{thm}
\begin{rmk}
Assumption (8.2) of \citet{gvdv} is omitted in the theorem statement above.
That is because it is satisfied by the Hellinger distance $d_H$ whenever the underlying collection of models $\calP$ is identifiable \citep[Proposition D.8]{gvdv}.
We will work in the Hellinger metric and assume identifiability throughout; see \eqref{identifiability} below.
\end{rmk}
\begin{rmk}
The premise of Theorem \ref{thm:ghosaletal:2000} is that the observations $\bfz^{(k)}$ are i.i.d.\ draws from $P_{ \theta_0 }$, but the prior $\Pi$ can depend on $k$.
This fact will be essential in later sections, where we apply Theorem \ref{thm:ghosaletal:2000} to phylogenetic trees in settings where the number of leaves in the latent tree increases with the amount of observed data.
\end{rmk}
The main (and the strongest) condition of the theorem is the so-called prior mass condition \eqref{eq:prior}, which ensures that we put sufficient prior probability mass near the true density $p_{ \theta_0 }$.
Here the proximity to the truth is measured through Kullback--Leibler divergence and the second Kullback--Leibler variation.
The condition is usually satisfied when the prior has exponential tails.
The entropy condition \eqref{eq:entropy} is the other main determinant of the contraction rate.
It bounds the complexity of the sieves $\calP_k$ which contain the majority of the prior probability mass.
The latter is ensured by the remaining mass condition \eqref{eq:remaining}.

\subsection{A state space of trees}\label{subsection:tree-state-space}

Consider a tree with a fixed number $n \in \bbN$ of leaves, and let $\TTT_n$ denote a finite set of possible topologies.
For each $T \in \TTT_n$, let $m(T) \in \bbN$ be the number of degrees of freedom required to specify the branch lengths of $T$, and let
\begin{equation*}
( \psi_T( x_1, \ldots, x_{ m( T ) } )_{ ( u, v ) } )_{ ( u, v ) \in T } =: ( \psi_T( \bfx_{ 1 : m( T ) } ) _{ ( u, v ) } )_{ ( u, v ) \in T }
\end{equation*}
be the mapping from variables specifying those degrees of freedom, $\bfx_{ 1 : m( T ) } \in ( 0, \infty )^{ m( T ) }$, to the lengths of branches $( ( u, v ) )_{ ( u, v ) \in T }$, where $( u, v )$ denotes a branch between two nodes $u$,  $v$ of the tree $T$.
When the length of the vector $\bfx_{ 1 : m( T ) }$ is clear from context, we will often omit the subscript and write $\bfx$ to lighten notation.
We will also denote by $I( T )$ the set of internal or non-leaf nodes of $T$.
For a given topology $T$, a branch $( u, v ) \in T$ will be called external if either $u$ or $v$ is a leaf, and the set of external branches will be denoted by $Ex(T)$.
Finally, let $m_n^* := \max_{ T \in \TTT_n }\{ m( T ) \}$.

\begin{figure}
	\centering
    \includegraphics[width=\textwidth]{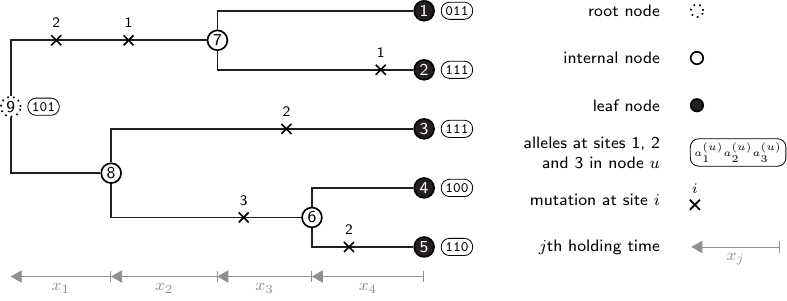}
    \caption{
		A binary, ultrametric tree $ T$ with $ n = 5 $ leaves.
		Time runs forwards from the root $ u_0 = 9 $ towards the leaves.
		We require $m(T) = 4$ holding times, $ x_1, \dotsc, x_4 $, to specify the branch lengths; for example, the branch from node 8 to 6 has length $ \psi_T(\bfx )_{(8,6)} = x_2 + x_3 $.
		The alleles at $ k = 3 $ sites taking values in $ H = \{0, 1\} $ are drawn independently from distribution $ r(\cdot) $ at the root and mutate according to a Markov jump process with rate matrix $ \bfQ $ along the branches of $ T $.
		We only observe an allele at the leaves so integrate over possible states at the internal nodes $ I(T) = \{6, \dotsc, 9\} $ in calculating the likelihood~\eqref{ultrametric_site_likelihood}.
	}
    \label{fig:phylogeny}
\end{figure}

This somewhat abstract description encompasses a broad range of tree models, such as those where the degree of tree nodes is not fixed and the number of edges varies by topology.
It also encompasses both rooted and unrooted trees, and can easily incorporate constraints, such as a requirement that trees be ultrametric.
\figref{fig:phylogeny} illustrates our notation on an ultrametric tree, and examples of the state space construction are presented in Sections \ref{subsection:kingman} and \ref{subsection:unconstrained}.
These examples consist of binary trees, but multifurcating trees could also be considered by taking $\TTT_n$ to be, say, the set of all unrooted trees with $n$ leaves.
Then $m( T )$ can be taken to be the number of edges in $T \in \TTT_n$, which depends on the number and size of multifurcations.
Similarly, $\psi_{ T }$ can be taken to be the identity function by imposing and ordering of edges $(u, v) \in T$ matching that of the vector $\bfx$.
In this way, the model selection problem of whether a tree contains multifurcations, as studied by \citet{zhangetal:2021}, can be recast as a single Bayesian inference problem with a prior distribution supporting both binary and multifurcating trees.

We assume the set $\TTT_n$, as well as the maps $m( \cdot )$ and $\psi_T( \cdot )$, are known a priori, and that $\psi_T$ is invertible for each $T$.
Then, our state space of trees is
\begin{equation*}
\Theta^{(n)} := \{ ( T, \bfx ) : T \in \TTT_n, \bfx \in ( 0, \infty )^{ m( T ) } \} .
\end{equation*}

\subsection{A model of evolution on trees}

Consider a fixed tree $( T, \bfx ) \in \Theta^{(n)}$.
We assume that the state, or allele of each leaf is observed at $k \in \bbN$ sites.
All sites evolve conditionally independently along the branches of the common genealogical tree.
At a given site, the root of the tree is assigned a random allele $a$ among a finite set of alleles $H$ from a distribution $r( \cdot )$.
Then, alleles are propagated towards leaves along branches according to the law of a Markov jump process with irreducible rate matrix $ \bfQ = (Q_{\alpha, \beta} : \alpha, \beta \in H) $.
These dynamics give rise to the single-site likelihood
\begin{equation}\label{ultrametric_site_likelihood}
P_{ ( T, \bfx ) }( \bfa ) := \sum_{ \bfa^* \in H^{ | I( T ) | } } r( a^{ ( u_0 ) } ) \prod_{ ( u, v ) \in T } \exp( \psi_T( \bfx )_{ ( u, v ) } \bfQ )_{ a^{ ( u ) } a^{ ( v ) } },
\end{equation}
where $\psi_T( \bfx )_{ ( u, v ) }$ is the length of the branch connecting nodes $u$ and $v$, $\bfa$ is the vector of observed alleles at the $n$ leaves, and the sum is over all non-leaf alleles in the extension, denoted $\bfa^*$, which tracks alleles at the internal nodes of the tree as well.
The allele at node $u$ is denoted by $a^{ ( u ) }$, with $u_0$ corresponding to the root.
\figref{fig:phylogeny} illustrates this mutation process on a rooted tree.

The likelihood \eqref{ultrametric_site_likelihood} depends on the position of the root, but can be extended to unrooted trees provided that the Markov jump process with rate matrix $\bfQ$ is $r$-reversible.
In that case, \cite{felsenstein:1981} showed that the position of the root does not affect the value of \eqref{ultrametric_site_likelihood}, and hence can be set arbitrarily.
He also presented an efficient algorithm for evaluating \eqref{ultrametric_site_likelihood}, which applies to both the reversible and rooted, non-reversible cases.
Likelihoods for multiple sites are products of terms of the form \eqref{ultrametric_site_likelihood}:
\begin{equation}\label{multisite_likelihood}
P_{ ( T, \bfx ) }( \bfa_1, \ldots, \bfa_k ) := \prod_{ i = 1 }^k \sum_{ \bfa_i^* \in H^{ | I( T ) | } } r( a_i^{ ( u_0 ) } ) \prod_{ ( u, v ) \in T } \exp( \psi_T( \bfx )_{ ( u, v ) } \bfQ )_{ a_i^{ ( u ) } a_i^{ ( v ) } }.
\end{equation}

In order to use Theorem \ref{thm:ghosaletal:2000}, we need upper bounds on the Kullback--Leibler divergence and centered second Kullback--Leibler variation corresponding to \eqref{multisite_likelihood}.
To that end, we show that both distances can be bounded by a Euclidean distance between $\bfx$ and $\bfx_0$ under suitable conditions on the lengths of tree branches, when the topology is fixed.
We first introduce some new notation.
Let $| \bfQ |$ be the matrix whose entries are the absolute values of those of $\bfQ$.
Let
\begin{align*}
\eta &:= \min_{ \alpha, \beta \in H }\{ Q_{ \alpha, \beta } : Q_{ \alpha, \beta } > 0 \}, \\
\gamma &:= \max_{ \alpha, \beta \in H }\{ | Q_{ \alpha, \beta } | \}, \\
w &:= \min\{ j \in \bbN : \min_{ \alpha, \beta } \{ ( | \bfQ |^j )_{ \alpha, \beta } \} > 0 \},
\end{align*}
be the least positive entry of $\bfQ$, the entry of $\bfQ$ with maximal magnitude, and the least number of steps in which the mutation Markov chain with generator $\bfQ$ can reach any state from any other with positive probability, respectively.
The last is guaranteed to exist since $\bfQ$ is a finite, irreducible matrix \citep[Chapter 1, Theorem 1.2.1]{norris:1997}.
For $f > 0$, $g > 0$, and a tree topology $T \in \TTT_n$, let
\begin{align*}
L_f( T ) &:= \big\{ \bfx \in ( 0, \infty )^{ m( T ) } : \min_{ ( u, v ) \in Ex( T ) }\{ \psi_T( \bfx )_{ ( u, v ) } \} > f \big\}, \\
U_g( T ) &:= \big\{ \bfx \in ( 0, \infty )^{ m( T ) } : \max_{ ( u, v ) \in T }\{ \psi_T( \bfx )_{ ( u, v ) } \} < g \big\},
\end{align*}
be the respective events that lengths of external branches are bounded below by $f$, and that the lengths of all branches are bounded above by $g$.
The complements of these events are denoted as $L_f^c(T)$ and $U_g^c(T)$, respectively.

\begin{lem}
\label{lem:distances}
Suppose $\psi_T( \cdot )$ is differentiable for each $T \in \TTT_n$, and that
\begin{equation*}
\max_{ T \in \TTT_n } \max_{ i \in \{ 1, \ldots, m( T ) \} } \Big| \frac{ \partial \psi_T }{ \partial x_i } ( \bfx ) \Big| =: \| \partial \psi^{(n)} \|_{ \infty } < \infty.
\end{equation*}
Define
\begin{equation}\label{Bn_definition}
B_n := A \| \partial \psi^{(n)} \|_{ \infty } ( n - 1 ) B^n := \frac{ 4 \gamma w! }{ \min_{ \alpha \in  H }\{ r( \alpha ) \} } \frac{ \| \partial \psi^{(n)} \|_{ \infty } ( n - 1 ) | H |^{ n - 1 } }{ \min_{ \alpha \in H, t \geq 0 } \{ \exp( \bfQ t )_{ \alpha \alpha } \}^{ n - 2 } \eta^{ n w } },
\end{equation}
where $A > 0$ and $B > 0$ are constants independent of $n$.
Then, for $0 < f < 1 / (2\eta),$
\begin{align}
\label{eq:klbound}
\KL( T_0, \bfx_0; T_0, \bfx ) &\leq B_n f ^{ - n w }  \| \bfx - \bfx_0 \|_2, \\
\label{eq:varbound}
V_0( T_0, \bfx_0; T_0, \bfx ) &\leq B_n^2 f ^{ - 2 n w }  \| \bfx - \bfx_0 \|^2_2,
\end{align}
as long as $\bfx,\bfx_0 \in L_f( T_0 )$.
\end{lem}

\section{A contraction theorem for posteriors on trees}\label{section:results}

Consider a prior distribution $\Pi^{(n)}$ on $\Theta^{(n)}$, and $n$ observed sequences at $k$ sites with respective allele configurations $( \bfa_1, \ldots, \bfa_k )$. From now on we assume the identifiability condition
\begin{equation}\label{identifiability}
P_{ \theta_0 } = P_{ \theta } \Rightarrow \theta_0 = \theta.
\end{equation}
\begin{rmk}
The identifiability condition \eqref{identifiability} is known to hold for a number of standard phylogenetic tree models \citep{steel/szekely:2007, steel/szekely:2009}.
A typical failure mode of  \eqref{identifiability} arises when entries of the matrix $\bfQ$ are part of the inference problem, whereupon a model only distinguishes compound products $\bfQ (\psi_T( \bfx ) )_{ ( u, v ) }$ of rates and branch lengths.
The typical solution is to reparametrize the model with these products directly.
\end{rmk}

Let the data-generating tree be $( T_0, \bfx_0 )$, and let
\begin{equation*}
\tilde f_n := \min_{ ( u, v ) \in Ex( T_0 ) }\{ \psi_{ T_0 }( \bfx_0 )_{ ( u, v ) }/2 \} \wedge 1 / (2\eta).
\end{equation*}
\begin{thm}
\label{thm:general_rooted}
Suppose \eqref{identifiability} holds, and that $n$ is either bounded as $k \to \infty$, or that if $n = n_k$ with $\lim_{ k \to \infty } n_k = \infty$, then there exists an infinite-leaf tree from which the sequence of data-generating trees are obtained as nested restrictions to $n_k$ leaves, for each $k$.
Suppose further that there exist sequences $\eps_k \to 0$ and $\zeta_k \to 0$ such that $k\zeta_k^2\to \infty$ and $\zeta_k\leq\eps_k$, sequences $f_k \to 0$  and $g_k \to \infty$ such that
\begin{equation}
\label{eq:entropy1}
\Big( \frac{ g_k }{ \eps_k^2 f_k^{ n w } } \Big)^{ m_n^* } \exp( - k \eps_k^2 ) \leq \frac{ 1 }{ | \TTT_n | ( 4 B_n )^{ m_n^* } },
\end{equation}
and a constant $C>0$ such that
\begin{align}
\label{eq:remaining1}
\Pi^{(n)}( \cup_{ T \in \TTT_n }\{ L_{ f_k }^c( T ) \cup U^c_{ g_k }( T ) \} ) \leq \exp( - ( C + 4 ) k \zeta_k^2 ),&\\
\label{eq:prior1}
\Pi^{(n)}( T_0 ) \Pi^{(n)}( \bfx \in L_{ \tilde f_n }( T_0 ) :
\|\bfx - \bfx_0 \|_2 \leq B_n^{-1}{\tilde f_n^{nw}}\zeta_k^2 \given T_0 ) \geq \exp( - C k \zeta_k^2 ).&
\end{align}
Then the posterior $\Pi^{(n)}( \cdot \given \bfa_1, \ldots, \bfa_k )$ converges to $( T_0, \bfx_0 )$ at rate $\eps_k$ as $k \to \infty$.
\end{thm}
Theorem \ref{thm:general_rooted} gives sufficient conditions for posterior consistency and a general tool for deriving posterior contraction rates. Notice that conditions \eqref{eq:entropy1}--\eqref{eq:prior1} can be matched one-to-one to the conditions in Theorem \ref{thm:ghosaletal:2000} by defining the sieves
\begin{equation}\label{sieves}
\Theta_k^{(n)} := \{ ( T, \bfx ) : T \in \TTT_n, \bfx \in  L_{ f_k }( T ) \cap U_{ g_k }( T ) \},
\end{equation}
which corresponds to selecting the trees that have a lower bound on the lengths of the external branches and an upper bound on all the branch lengths. Then \eqref{eq:entropy1} ensures that the entropy condition \eqref{eq:entropy} is satisfied, while \eqref{eq:remaining1} and \eqref{eq:prior1} correspond to the remaining and the prior mass conditions respectively.

The requirement of a nested sequence of data-generating trees in Theorem \ref{thm:general_rooted} ensures that there is a well-defined limit object around which the posterior can converge as $n, k \to \infty$.
In order to prove convergence, we require $k$ to grow rapidly enough with $n$ so that \eqref{eq:entropy1}--\eqref{eq:prior1} can be shown to hold.
The examples in Section \ref{section:examples} illustrate this requirement.

\section{Examples}\label{section:examples}

In this section we show that the conditions of Theorem \ref{thm:general_rooted} hold for two popular priors: the Kingman coalescent prior for ultrametric trees, and a uniform topology prior with independent exponential branch lengths on unconstrained trees.
Both settings will utilize the following elementary lemma.
\begin{lem}
\label{lem:branches_bounds}
Let $( X_1, \ldots, X_n )$ be random variables with marginal laws $X_i \sim \operatorname{Exp}( \lambda_i )$.
For any $f,g > 0$,
\begin{align*}
\bbP( X_i < f) &\leq \lambda_i f \Bigg[ 1 + \frac{\lambda_i}{2} f \Bigg], \\
\bbP(\max\{ X_1, \ldots, X_n \}  > g ) &\leq n \exp( - \min\{ \lambda_1, \ldots, \lambda_n \} g ).
\end{align*}
\end{lem}

\subsection{Kingman coalescent prior}\label{subsection:kingman}

Let $\TTT_n^K$ be the set of \emph{ranked topologies} of rooted binary trees with $n$ labeled leaves, of which there are $\binom{ n }{ 2 } \binom{ n - 1 }{ 2 } \cdots \binom{ 2 }{ 2 } = 2^{ - n + 1 }( n! )^2 / n \sim 4 \pi [ n / ( \sqrt{2} e ) ]^{2 n}$.
A ranked topology specifies both the binary mergers which take place, and the order in which they occur.
To ensure the prior places full mass on ultrametric trees, it is convenient to parametrize branch lengths by the $n - 1$ holding times between mergers, $\bfx \in ( 0, \infty )^{ n - 1 }$, so $m( T ) \equiv n - 1$ for all $ T $, and
\begin{equation}\label{eq:kingman_psi}
\psi_T( \bfx )_{ ( u, v ) } = \sum_{ j \in \{ 1, \ldots, n - 1 \} : x_j \in ( u, v ) } x_j,
\end{equation}
where $x_j \in ( u, v )$ means that holding time $x_j$ contributes to the length of branch $( u, v )$.
See \figref{fig:phylogeny} for an illustration.
Let
\begin{equation*}
\Pi^{ K, n }( T, \bfx ) := \exp\Bigg( - \sum_{ i = 1 }^{ n - 1 } \binom{ n - i + 1}{2} x_i \Bigg)
\end{equation*}
denote the Kingman coalescent prior \citep{Kingman82a} on $\Theta^{ K, n } := \TTT_n^K \times ( 0, \infty )^{ n - 1 }$.
It is obtained by initializing $n$ disconnected leaf lineages from time zero and merging each pair at rate 1 until their most recent common ancestor is reached.
The resulting distribution on trees is uniform over ranked topologies, and has independent holding times $X_i \sim \operatorname{Exp}( \binom{ i + 1 }{ 2 } )$ for $i \in \{ 1, \ldots, n - 1 \}$ governing its branch lengths.

Corollary \ref{thm:main} below considers a sequence of trees tending to an infinite-leaf limit.
The fact that such a limit exists is a well-known consistency property of the Kingman coalescent \citep[Theorem 3]{Kingman82a}, and we denote its law by $\Pi^K$.
Indeed, a Kingman coalescent tree with $n$ leaves can be extended to one with $n + 1$ leaves by splitting a randomly chosen leaf into two, and extending all external branches (including the two zero-length ones just created by the split) by an independent $X_n \sim \text{Exp}\big( \binom{n+1}{2} \big)$.
This ordering of holding times, where $x_{n-1}$ is the time until the first merger when there are $n$ leaves, and $x_1$ is the holding time between the last two mergers, is depicted in Figure \ref{fig:phylogeny}.
For an infinite-leaf tree $( T, \bfx ) \in \TTT_{ \bbN }^K \times ( 0, \infty )^{ \bbN }$, let $( T, \bfx )|_n = ( T|_n, \bfx|_n )$ be the restriction of the ranked topology to the first $n$ leaves with branch lengths determined by holding times $\bfx|_n = ( x_1, \ldots, x_{n - 1} )$.

\begin{cor}
\label{thm:main}
Suppose $( T_0, \bfx_0 )$ is an infinite-leaf tree in the support of $\Pi^K$, and that the data-generating sequence of trees is given by the restrictions $( T_0, \bfx_0 )|_n$.
Let the observed sites $( \bfa_1, \ldots \bfa_k )$ be coupled such that for any $n$, trees $( T_0, \bfx_0 )|_n$ and $( T_0, \bfx_0 )|_{ n + 1 }$ use the same Poisson process on $( T_0, \bfx_0 )|_n$ and $( ( T_0, \bfx_0 )|_{ n + 1 } ) |_n$, and a further independent Poisson process on the difference $( T_0, \bfx_0 )|_{ n + 1 } \setminus ( T_0, \bfx_0 )|_n$.
Suppose $n \equiv n_k$ increases slowly enough that for any sufficiently large $k$, and a constant $C > 0$ independent of $k$ and $n$, we have $\log k \geq n^2 / C^2$.
Then posterior distributions $\Pi^{ K, n }( \cdot \given \bfa_1, \ldots, \bfa_k )$ obtained from priors $\Pi^{ K, n }$ concentrate around $( T_0, \bfx_0 ) $ at rate $\eps_k\asymp k^{-1/2}(\log k)^{3 / 4 + \beta / 2 }$ for any $\beta > 1$ as $k \to \infty$.
\end{cor}

\begin{figure}
	\centering
	\begin{subfigure}[t]{\textwidth}
		\centering
	    \includegraphics[width=\textwidth]{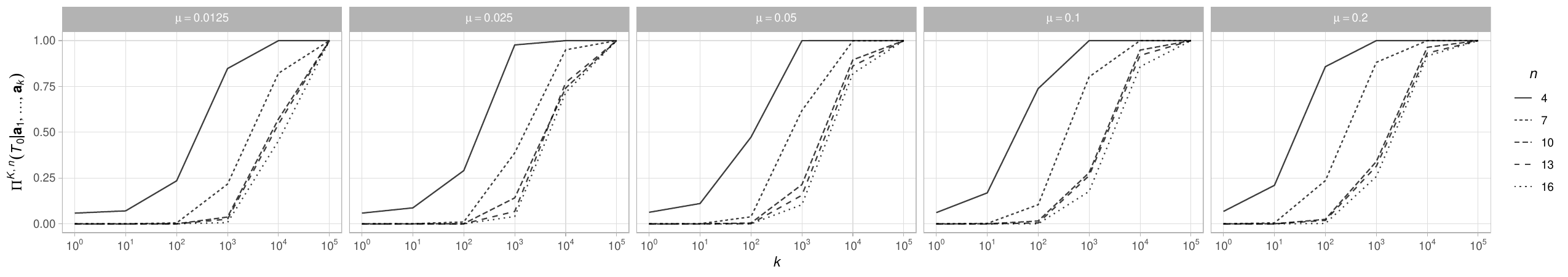}
	    \caption{
            Posterior support for the true tree topology \( T_0 |_ n \) averaged across 100 replicate data sets for each \( n \), \( \mu \) and \( k \).
			A figure depicting the posterior support from each data set is contained in the appendix.
		}
	    \label{fig:kingman-support}
	\end{subfigure}
	\begin{subfigure}[t]{\textwidth}
		\centering
	    \includegraphics[width=\textwidth]{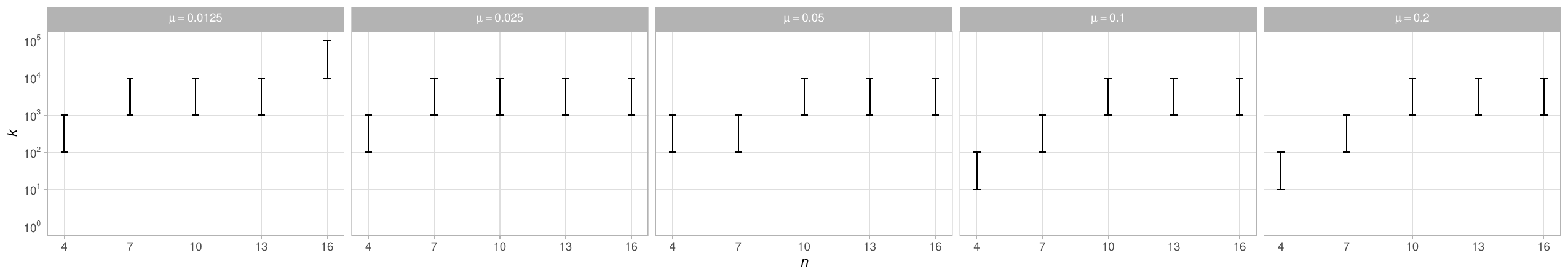}
	    \caption{
            Intervals for $ k $ on which the curves in \figref{fig:kingman-support} first exceed 0.5.
		}
	    \label{fig:kingman-threshold}
	\end{subfigure}
	\caption{
		Posterior support on the true rooted tree topology $ T_0 |_n $ under a Kingman coalescent prior $ \Pi^{K, n} $ as the number of leaves $ n $, number of sites $ k $ and mutation rate $ \mu $ vary.
        The sequence of trees $ (T_0, \bfx_0) |_n \sim \Pi^{K, n} $ was generated according to the algorithm in Section~\ref{subsection:kingman} and binary data sets were sampled according to a Markov model with mutation rate \( \mu \).
        See Section~\ref{subsection:further-details-experiments} for details of the inference.
	}
    \label{fig:kingman}
\end{figure}

To illustrate Corollary~\ref{thm:main}, \figref{fig:kingman-support} depicts the posterior support for the true tree $ T_0 $ as the number of leaves $ n $ and sites $ k $ increase in simple experiments with binary alleles and symmetric mutation rates \( \mu = Q_{0, 1} = Q_{1, 0} \). \figref{fig:kingman-threshold} displays the intervals on which the curves in \figref{fig:kingman-support} first exceed 0.5.

\subsection{Uniform prior}\label{subsection:unconstrained}

Let $\TTT_n^U$ be the set of unrooted, binary trees with $n$ labeled leaves.
There are
\begin{equation}\label{uniform_topologies}
| \TTT_n^U | = (2 n - 5 )!! := 2^{ 3 - n } ( 2 n - 5 )! / ( n - 3 )! \sim 2^{ 3 - n } \sqrt{ \frac{ 2 n - 5 }{ n - 3 } } \frac{ ( 2 n - 5 )^{ 2 n - 5 } }{ e^{ n - 8 } ( n - 3 )^{ n - 3 } }
\end{equation}
labeled topologies for such trees \citep{felsenstein:1978}, each requiring the specification of $2 n - 3$ branch lengths.
Thus, we set $m( T ) \equiv 2 n - 3$, take the variables $( x_1, \ldots, x_{ 2 n - 3 } ) \in ( 0, \infty )^{ 2 n - 3 }$ to parametrize branch lengths directly, and set $\psi_T( \bfx ) = \bfx$.

For $\lambda > 0$, let
\begin{equation*}
\Pi^{ U, n }( T, \bfx ) := \frac{ \lambda^{ 2 n - 3 } }{ | \TTT_n^U | }  \exp\Bigg( - \lambda \sum_{ i = 1 }^{ 2 n - 3 } x_i \Bigg)
\end{equation*}
be the prior on $\Theta^{ U, n } := \TTT_n^U \times ( 0, \infty )^{ 2 n - 3 }$ which is uniform over topologies and assigns independent exponentially distributed lengths with rate $\lambda$ to branches.
This prior was used, for example, by \cite{whidden/matsen:2015} to benchmark MCMC algorithms on phylogenetic trees.

Corollary~\ref{thm:unconstrained_main} below considers a sequence of trees tending to an infinite-leaf limit.
A consistent sequence of finite restrictions $( T, \bfx )|_n$ of an infinite-leaf tree can be constructed via the \emph{sequential construction} of \cite{aldous:1993}:
\begin{enumerate}
\item To extend $( T, \bfx ) \in \Theta^{U, n}$ to $\Theta^{U, n + 1}$, choose uniformly at random a branch $b$ on $T$.
Attach two new branches at a uniformly chosen end of $b$: one of length $x_{2 n - 2} \sim \text{Exp}(\lambda)$ connecting $b$ to its previous endpoint, and one of independent length $x_{2 n - 1} \sim \text{Exp}(\lambda)$ connecting the resulting node to a new leaf.
\item To restrict $( T, \bfx ) \in \Theta^{U, n + 1}$ to $\Theta^{U, n}$, choose a leaf uniformly at random and remove it.
Remove the resulting degree-2 internal node by deleting the adjacent edge whose length has a higher index in $\bfx$ and attaching its end points to each other directly.
\end{enumerate}
Let $\Pi^U$ be the law of the infinite-leaf projective limit obtained from the above procedure, which exists by the Kolmogorov extension theorem.
For an infinite-leaf tree $( T, \bfx ) \in \TTT_{ \bbN }^U \times ( 0, \infty )^{ \bbN }$, let $( T, \bfx )|_n = ( T|_n, \bfx|_n )$ be the restriction of the unrooted topology to the first $n$ leaves with branch lengths $\bfx|_n = ( x_1, \ldots, x_{ 2 n - 3 } )$.

\begin{cor}\label{thm:unconstrained_main}
Suppose $( T_0, \bfx_0 )$ be an infinite-leaf tree in the support of $\Pi^U$, and that the data-generating sequence of trees is given by the restrictions $( T_0, \bfx_0 )|_n$.
Let the observed sites $( \bfa_1, \ldots \bfa_k )$ be coupled such that for any $n$, trees $( T_0, \bfx_0 )|_n$ and $( T_0, \bfx_0 )|_{ n + 1 }$ use the same Poisson process on $( T_0, \bfx_0 )|_n$ and $( ( T_0, \bfx_0 )|_{ n + 1 } ) |_n$, and a further independent Poisson process on the difference $( T_0, \bfx_0 )|_{ n + 1 } \setminus ( T_0, \bfx_0 )|_n$.
Suppose $n \equiv n_k$ increases slowly enough that for any sufficiently large $k$, and a constant $C > 0$ independent of $k$ and $n$, we have $\log k \geq n^2 / C^2$.
Then posterior distributions $\Pi^{ U, n }( \cdot \given \bfa_1, \ldots, \bfa_k )$ obtained from priors $\Pi^{ U, n }$ concentrate around $( T_0, \bfx_0 ) $ at rate $\eps_k\asymp k^{-1/2}(\log k)^{3 / 4 + \beta / 2 }$ for any $\beta > 1$ as $k \to \infty$, provided the rate matrix $\bfQ$ is $r$-reversible.
\end{cor}

\begin{figure}
	\centering
	\begin{subfigure}[t]{\textwidth}
		\centering
	    \includegraphics[width=\textwidth]{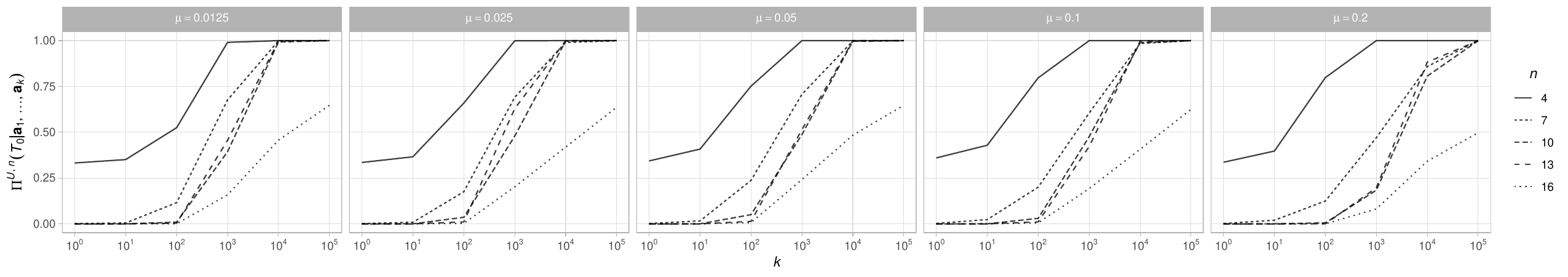}
	    \caption{
            Posterior support for the true tree topology \( T_0 |_ n \) averaged across 100 replicate data sets for each \( n \), \( \mu \) and \( k \).
			A figure depicting the posterior support from each data set is contained in the appendix.
		}
	    \label{fig:uniform-support}
	\end{subfigure}
	\begin{subfigure}[t]{\textwidth}
		\centering
	    \includegraphics[width=\textwidth]{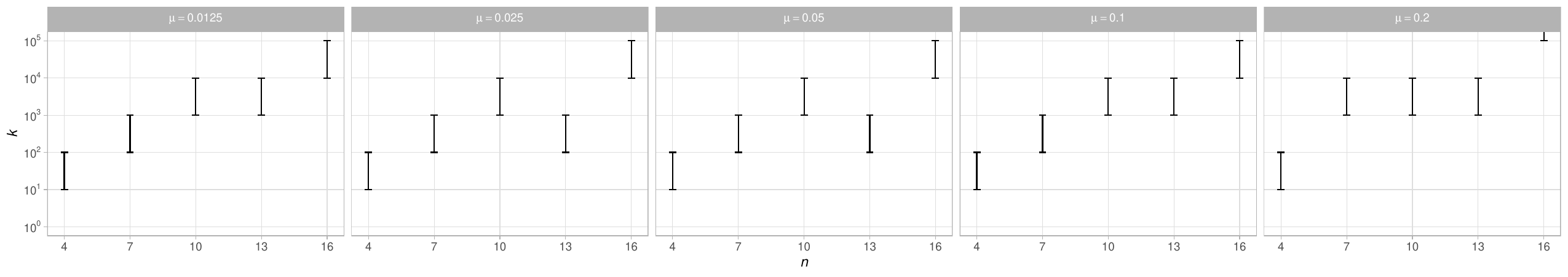}
	    \caption{
            Intervals for $ k $ on which the curves in \figref{fig:uniform-support} first exceed 0.5.
		}
	    \label{fig:uniform-threshold}
	\end{subfigure}
	\caption{
        Posterior support on the true unrooted tree topology $ T_0 |_n $ under a prior $ \Pi^{U, n} $ uniform across topologies and exponential on branch lengths as the number of leaves $ n $, number of sites $ k $ and mutation rate $ \mu $ vary.
        A sequence of trees $ (T_0, \bfx_0) |_n \sim \Pi^{U, n} $ was generated according to the algorithm in Section~\ref{subsection:unconstrained}
		and binary data sets were sampled according to a Markov model with mutation rate \( \mu \).
        See Section~\ref{subsection:further-details-experiments} for details of the inference.
	}
    \label{fig:uniform}
\end{figure}

To illustrate Corollary~\ref{thm:unconstrained_main}, \figref{fig:uniform-support} shows the posterior distributions on unrooted tree topologies concentrating in a simple experiment with binary alleles and symmetric mutation rates \( \mu = Q_{0, 1} = Q_{1, 0} \).
\figref{fig:uniform-threshold} displays the intervals for which the posterior support on the true topology first exceeds 0.5.

\subsection{Further details on experiments}
\label{subsection:further-details-experiments}

\texttt{RevBayes} \citep{hohna2016revbayes} was used to perform inference via Markov chain Monte Carlo.
Each chain was initialized randomly and the first $ 10^6 $ iterations discarded as warm up, then $ 10^3 $ samples were drawn at intervals of $ 10^3 $ iterations and used to estimate the posterior support on the true tree topology.
Mixing and convergence was assessed visually via trace plots of the log-likelihood at each iteration and computing the overall acceptance rate of proposals.
Plots were created using \texttt{ggplot2} \citep{wickham16} in R \citep{r2024}.
Code to recreate the experiments is available at \url{github.com/lukejkelly/TreeConsistency}.
The routines to generate trees and allele data make use of \texttt{ape} \citep{paradis2019ape} and \texttt{phangorn} \citep{schliep2011phangorn}.
The appendix contains additional figures and details of the experiments.

\section{Discussion}\label{section:discussion}

We have demonstrated the consistency of Bayesian tree reconstruction procedures in a wide range of settings, and under mild technical conditions.
In particular, we do not require any discretisation or boundedness assumptions on branch lengths, in contrast with earlier results.
Our results apply to rooted and unrooted trees, unconstrained and ultrametric (or otherwise constrained) trees, and to both binary and more flexible tree models.

We presented several prototypical examples of tree reconstruction problems, and obtained explicit convergence rates for them which match known frequentist rates, obtained under stronger assumptions.

While our results apply to the setting in which the number of leaves, and hence the size of the tree, diverge with the number of replicate sites, we require the number of sites to be super-exponential in the number of leaves.
Under stronger assumptions on branch lengths (and sometimes the mutation model), a polynomial number of sites suffices for accurate maximum likelihood reconstruction \citep{roch/sly:2017, dinhetal:2018}.
Since maximum likelihood reconstruction corresponds to maximum a posteriori estimation under a uniform prior when branch lengths are bounded, it is clear that similar refinements can be obtained for particular models in the Bayesian setting.
Our numerical results in Figures \ref{fig:kingman-threshold} and \ref{fig:uniform-threshold} are consistent with sub-exponential growth of $k$ in $n$, strongly suggesting that the super-exponential growth rate which we require for our method of proof can be improved, at least in particular cases.

Our rates pertain to convergence on the space of distributions, $P_{ \theta } \to P_{ \theta_0 }$.
In addition, it is also of interest to obtain rates for convergence of the corresponding parameter estimation problem, $\theta \to \theta_0$.
In the context of genetics, the latter corresponds to learning the ancestral tree, while the former corresponds to learning the distribution of alleles at the $n$ leaves of the tree at an arbitrary site.
The identifiability assumption \eqref{identifiability} guarantees that $P_{ \theta } \to P_{ \theta_0 } \Rightarrow \theta \to \theta_0$, but the rate of the latter depends on the fluctuations of the $\theta \mapsto P_{ \theta }$ mapping and can, in principle, be arbitrarily slow.
When they are available, model-specific results bounding $\| \theta - \theta_0 \|$ by $\| P_{ \theta } - P_{ \theta_0 } \|$ in suitable norms can be used to obtain rates for the $\theta \to \theta_0$ convergence as well; see \citet[Eq.\ (3)]{zhangetal:2021} for an example.

\section{Proofs of main results}\label{section:proofs}

\subsection{Proof of Theorem \ref{thm:general_rooted}}
We are going to match conditions \eqref{eq:entropy1}--\eqref{eq:prior1} to the conditions \eqref{eq:remaining}--\eqref{eq:prior} in Theorem \ref{thm:ghosaletal:2000}.
Let the subsets $\Theta_k^{ ( n ) }\subset \Theta^{ ( n ) }$ be given by \eqref{sieves}, and define $\calP_k^{ ( n ) }=\{P_\theta : \theta\in\Theta_k^{ ( n ) }\}$.
Then the remaining mass condition \eqref{eq:remaining} is clearly satisfied when \eqref{eq:remaining1} holds.

For the entropy condition \eqref{eq:entropy}, recall that $d^2_H( p, q )\leq \KL( p ; q )$ for any pair of probability densities $p, q$ \cite[Lemma B.1(iii)]{gvdv}.
Therefore, we can use \eqref{eq:klbound} and the identifiability condition \eqref{identifiability} to bound the covering number as
\begin{align*}
N(\eps_k/2,\calP_k^{ ( n ) }, d_H) \leq N(\eps^2_k/4,\calP_k^{ ( n ) }, \KL) &\leq |\TTT_n| N\left(\frac{ \eps^2_k f_k^{ n w } }{ 4 B_n }, [ 0, g_k ]^{ m_n^* }, \| \cdot \|_2 \right) \\
&\leq | \TTT_n | \left( \frac{ 4 B_n g_k }{ \eps_k^2 f_k^{ n w } } \right)^{ m_n^* }.
\end{align*}
Therefore, for \eqref{eq:entropy} to hold we require the following, which is equivalent to \eqref{eq:entropy1}:
\[
\frac{ g_k^{ m_n^* } }{ \eps_k^{ 2 m_n^* } f_k^{ n w m_n^* } } \leq \frac{ \exp( k \eps_k^2 ) }{ | \TTT_n | ( 4 B_n )^{ m_n^* } }.
\]
It remains to obtain \eqref{eq:prior}. Let $\theta_0=(T_0,\bfx_0)$ be the true parameter.
Then for $\bfx \in L_{ \tilde f_n }( T_0 )$ we can use  \eqref{eq:klbound} and \eqref{eq:varbound} to obtain
\begin{align*}
&\Pi^{(n)}( ( T, \bfx ) \in \Theta : \KL( P_{\theta_0}; P_\theta )\leq \zeta_k^2,\, V_0( P_{\theta_0}; P_\theta )\leq \zeta_k^2 ) \\
&\geq \Pi^{(n)}( T_0 ) \Pi^{(n)}( \bfx \in ( 0, \infty )^{ m( T_0 ) } : \KL( P_{(T_0,\bfx_0)}; P_{(T_0,\bfx)} ) \leq \zeta_k^2,\, V_0( P_{(T_0,\bfx_0)}; P_{(T_0,\bfx)} )\leq \zeta_k^2 \given T_0 )\\
&\geq \Pi^{(n)}( T_0 ) \Pi^{(n)}\Bigg( \bfx \in L_{ \tilde f_n }( T_0 ): \frac{ B_n }{ \tilde f_n^{nw} } \| \bfx - \bfx_0 \|_2 \leq \zeta_k^2, \frac{ B_n^2 }{ \tilde f_n^{2nw} } \| \bfx - \bfx_0 \|^2_2 \leq \zeta_k^2 \given T_0 \Bigg).
\end{align*}
For $k$ large enough that $\zeta_k\leq 1$, the condition $B_n^2\| \bfx - \bfx_0\|^2_2\leq \tilde f_n^{2nw}\zeta_k^2$ is automatically satisfied once $B_n \|\bfx-\bfx_0\|_2\leq \tilde f_n^{nw}\zeta_k^2$.
Therefore,
\begin{align*}
&\Pi^{(n)}( ( T, \bfx ) \in \Theta : \KL( P_{\theta_0}; P_\theta )\leq \zeta_k^2,\, V_0( P_{\theta_0}; P_\theta )\leq \zeta_k^2 ) \\
&\geq \Pi^{(n)}( T_0 ) \Pi^{(n)}( \bfx \in L_{ \tilde f_n }( T_0 ) : \|\bfx - \bfx_0 \|_2 \leq B_n^{-1}{\tilde f_n^{nw}}\zeta_k^2 \given T_0 ),
\end{align*}
at which point \eqref{eq:prior1} implies \eqref{eq:prior}.

\subsection{Proof of Corollary \ref{thm:main}}

Before proving Corollary \ref{thm:main}, we collect some useful properties of the infinite-leaf coalescent tree $( T, \bfx ) \sim \Pi^K$ in a preparatory lemma.

\begin{lem}\label{lem:kingman_shape}
With $\Pi^K$-probability one, $\| \bfx \|_1 < \infty$, and only finitely many of the events $\{ x_{ n - 1 } < n^{ -4 } \}_{ n \geq 2 }$ happen.
\end{lem}
\begin{proof}
Recall that $( x_1, x_2, \ldots )$ are independent and $x_i \sim \text{Exp}\big( \binom{i+1}{2} \big)$.
For the first claim,
\begin{equation*}
\bbE[ \| \bfx \|_1 ]  = \bbE\Bigg[ \sum_{ i = 1 }^{ \infty } x_i \Bigg] = \sum_{ i = 1 }^{ \infty } \frac{ 2 }{ i ( i + 1 ) } = 2 \sum_{ i = 1 }^{ \infty } \Big( \frac{ 1 }{ i } - \frac{ 1 }{ i + 1 } \Big) = 2,
\end{equation*}
where the interchange of sum and expectation can be justified by Tonelli's theorem.
Hence, $\Pi^K( \| \bfx \|_1 = \infty ) = 0$.
For the second claim, Lemma \ref{lem:branches_bounds} yields
\begin{equation*}
\sum_{ n = 2 }^{ \infty }\Pi^K( x_{ n - 1 } < n^{ -4 } ) \leq \sum_{ n = 2 }^{ \infty } \frac{ n^2 }{ 2 } n^{ -4 } \Big( 1 + \frac{ n^2 }{ 4 } n^{ -4 } \Big) = \sum_{ n = 2 }^{ \infty } \frac{ 5 }{ 8 n^2 } < \infty.
\end{equation*}
The claim follows by the Borel--Cantelli lemma.
\end{proof}

For $\beta > 1$, $\delta > 0$, and the constant $C > 0$ from the statement of Corollary \ref{thm:main}, define
\begin{align*}
\zeta_k^2 &= \frac{ C ( \log k )^{ 3 / 2 } }{ k }, \\
\eps_k^2 &= \frac{ C ( \log k )^{ 3 / 2 + \beta } }{ k },\\
f_k &= \frac{ 1 }{ n( n - 1 ) } \exp( - ( C + 4 ) k \zeta_k^2 ),\\
g_k &= ( C + 4 + \delta ) k \zeta_k^2.
\end{align*}
Then, using $\log k \geq n^2 / C^2$ yields
\begin{align*}
&\frac{ g_k^{ n - 1 } }{ \eps_k^{ 2 ( n - 1 ) } f_k^{ n w ( n - 1 ) } } \exp( - k \eps_k^2 ) \\
&\leq  \frac{ ( C + 4 + \delta )^n n^{ 2 w n^2 } }{ ( \log k )^{ \beta ( n - 1 ) } } \exp( [ n + C ( C + 4 ) w n^2 ( \log k )^{ 1 / 2 } - C ( \log k )^{ 1 / 2 + \beta } ] \log k ) \\
&\leq \frac{ ( C + 4 + \delta )^n n^{ 2 w n^2 } }{ ( \log k )^{ \beta ( n - 1 ) } } \exp( [ n + \{ C ( C + 4 ) w n^2 - C ( \log k )^{ \beta } \} ( \log k )^{ 1 / 2 } ] \log k ) \\
&\leq \frac{ ( C + 4 + \delta )^n n^{ 2 w n^2 } }{ ( \log k )^{ \beta ( n - 1 ) } } \exp( [ n + \{ C ( C + 4 ) w n^2 - C^{ 1 - 2 \beta } n^{ 2 \beta }  \} ( \log k )^{ 1 / 2 } ] \log k ).
\end{align*}
Since $\beta > 1$, the factor inside the $\{ \}$-braces is negative for large enough $n$.
Hence, for such a large enough $n$,
\begin{align}
&\frac{ g_k^{ n - 1 } }{ \eps_k^{ 2 ( n - 1 ) } f_k^{ n w ( n - 1 ) } } \exp( - k \eps_k^2 )  \notag \\
&\leq \frac{ ( C + 4 + \delta )^n n^{ 2 w n^2 } }{ ( \log k )^{ \beta ( n - 1 ) } } \exp( [ n + \{ C ( C + 4 ) w n^2 - C^{ 1 - 2 \beta } n^{ 2 \beta } \} ( n / C ) ] \log k ) \notag \\
&\leq \frac{ ( C + 4 + \delta )^n }{ ( n^2 / C^2 )^{ \beta ( n - 1 ) } } \exp( 2 w n^2 \log n + [ n + \{ C ( C + 4 ) w n^2 - C^{ 1 - 2 \beta } n^{ 2 \beta } \} ( n / C ) ] \log k ) \notag \\
&\leq \frac{ ( C + 4 + \delta )^n }{ ( n^2 / C^2 )^{ \beta ( n - 1 ) } } \exp( [ 2 w \log n + n / C^2 + C^{ -2 } ( C + 4 ) w n^3 - C^{ -2 ( 1 + \beta ) } n^{ 1 + 2 \beta } ] n^2 ), \label{cor1_lhs}
\end{align}
where the last line follows because the negative sign dominates for large enough $n$ since $\beta > 1$.

Recalling the definition of $B_n$ in \eqref{Bn_definition} and the fact that $\| \partial \psi \|_{ \infty } \leq n - 1$ in \eqref{eq:kingman_psi}, the right-hand side of \eqref{eq:entropy1} is bounded below as
\begin{align}
\frac{ 1 }{ | \TTT_n^K | ( 4 B_n )^{ n - 1 } } &\sim \frac{ 1 }{ 4 \pi } \Bigg( \frac{ \sqrt{ 2 } e }{ n } \Bigg)^{ 2n } \frac{ 1 }{ A \| \partial \psi^{ ( n ) } \|_{ \infty } ( n - 1 ) B^n } \notag \\
&\geq \frac{ 1 }{ 4 \pi A n^2 } \Bigg( \frac{ 2 e^2 }{ n^2 B } \Bigg)^n. \label{cor1_rhs}
\end{align}
When $\beta > 1$, \eqref{cor1_lhs} decays faster than \eqref{cor1_rhs} as $n \to \infty$.
Hence \eqref{eq:entropy1} holds for any sufficiently large $n$.

In order to show \eqref{eq:remaining1} we use Lemma \ref{lem:branches_bounds} twice.
Firstly, the holding time $x_{ n - 1 } \sim \operatorname{Exp}\big( \binom{ n }{ 2 } \big)$ contributes to the lengths of all external branches in $( T, \bfx )|_n$.
Hence $\cup_{ T \in \TTT_n^K } L_f^c( T ) = \{ x_{ n - 1 } < f \}$.
Thus, by Lemma \ref{lem:branches_bounds} and the definition of $f_k$,
\[
\Pi^{ K, n }( \cup_{ T \in \TTT_n^K } L^c_{ f_k }( T ) )\leq \frac{1}{2} \exp( - ( C + 4 ) k \zeta_k^2 ) \Big( 1 + \frac{ 1 }{ 4 } \Big) =\frac{ 5 }{ 8 } \exp( - ( C + 4 ) k \zeta_k^2 ).
\]
On the other hand, $\cup_{ T \in \TTT_n^K } U_g^c( T ) \subseteq \{ \max_{ i \in \{ 1, \ldots, n - 1 \} }\{ x_i \} > g \}$.
Again by Lemma \ref{lem:branches_bounds},
\begin{align*}
\Pi^{ K, n }( \cup_{ T \in \TTT_n^K } U^c_{ g_k }( T ) ) \leq ( n - 1 ) \exp( - g_k ) &= ( n - 1 ) \exp( - ( C + 4 + \delta ) k \zeta_k^2 ) \\
&\leq \frac{ 3 }{ 8 } \exp( - ( C + 4 ) k \zeta_k^2 ),
\end{align*}
when $k$ is large enough that $- \delta C ( \log k )^{ 3 / 2 } + \log( n - 1 ) \leq \log ( 3 / 8 )$, which exists because $k$ grows superlinearly in $n$.
Therefore, \eqref{eq:remaining1} is satisfied.

The left-hand side of \eqref{eq:prior1} does not depend on $T_0$ because ranked topologies and branch lengths are independent under $\Pi^K$ and $\Pi^{K, n}$, and the marginal distribution over ranked topologies is uniform.
Let $B_n' := \tilde{f}_n^{nw} / (\sqrt{n}B_n)$.
We have that $\| \partial \psi^{ ( n ) } \|_{ \infty } \geq 1$, and by definition $\tilde{ f }_n \leq 1 / ( 2 \eta )$.
Hence
\begin{align}
B_n' \zeta_k^2 = \frac{ \tilde{ f }_n^{ n w } }{ \sqrt{ n } B_n } \frac{ C ( \log k )^{ 3 / 2 } }{ k } &\leq \frac{ \tilde{ f }_n^{ n w } }{ n^{ 3/2 } A \| \partial \psi^{(n)} \|_{ \infty } B^n } \frac{ n^3 }{ C^2 \exp( n^2 / C^2 ) } \notag\\
&\leq \frac{ n^{ 3/2 } }{ A C^2 B^n (2 \eta)^{ w m } \exp( n^2 / C^2 ) } \to 0. \label{Bn_decay}
\end{align}
Again by Lemma \ref{lem:kingman_shape}, $\min\{ x_1, \ldots, x_{ n - 1 } \} \to 0$ no faster than $n^{-4}$, which is slower than the decay in \eqref{Bn_decay}.
Hence we can take a large enough $n$ and $k$, so that $x_{ 0, i } - B_n' \zeta_k^2 > 0$ for each $i \in \{ 1, \ldots, n - 1 \}$, and $| x_{ n - 1 } - x_{ 0, n - 1 } | \leq B_n' \zeta_k^2 \Rightarrow x_{ n - 1 } \geq \tilde f_n$, whereupon
\begin{align*}
&\Pi^{ K, n }( \bfx : x_1\geq \tilde f, \, \|\bfx - \bfx_0|_n\|_2\leq B_n^{-1}\tilde f^{nw}\zeta_k^2 ) \\
&=\int_{\bfx \in \bbR_+^{ n - 1 } : \, x_1 \geq \tilde f, \| \bfx - \bfx_0|_n \|_2 \leq \sqrt{n}B_n'\zeta_k^2} \Bigg[ \prod_{ i = 1 }^{ n - 1 } \binom{ i }{ 2 } \exp\Bigg( - \binom{ i }{ 2 } x_i \Bigg) \Bigg] \rmd \bfx \\
&\geq \int_{\bfx \in \bbR_+^{ n - 1 } : \, | x_i - x_{ 0 , i } | \leq B_n' \zeta_k^2 } \Bigg[ \prod_{ i = 1 }^{ n - 1 } \binom{ i }{ 2 } \exp\Bigg( - \binom{ i }{ 2 } x_i \Bigg) \Bigg] \rmd \bfx \\
&=\prod_{ i = 1 }^{ n - 1 }\Bigg[ \exp\Bigg( - \binom{ i }{ 2 } ( x_{ 0 , i } - B_n' \zeta_k^2 ) \Bigg) - \exp\Bigg( - \binom{ i }{ 2 } ( x_{ 0 , i } + B_n' \zeta_k^2 ) \Bigg) \Bigg].
\end{align*}
We use a Taylor expansion to show that for any $i\in\{ 1 , \dots, n - 1 \},$
\[
\exp\Bigg( - \binom{ i }{ 2 } ( x_{ 0, i } - B_n' \zeta_k^2 ) \Bigg) - \exp\Bigg( - \binom{ i }{ 2 } ( x_{ 0, i } + B_n' \zeta_k^2 ) \Bigg) \geq 2 \binom{ i }{ 2 } B_n' \zeta_k^2 \exp\Bigg( - \binom{ i }{ 2 } x_{ 0 , i } \Bigg).
\]
Therefore,
\begin{align*}
\Pi^{ K, n }( \bfx \in ( 0, \infty )^{ n - 1 } : x_1 \geq \tilde f, &\|\bfx - \bfx_0|_n \|_2 \leq \sqrt{n} B_n' \zeta_k^2 ) \\
&\geq ( B_n' )^{ n - 1 } \frac{ ( n! )^2 }{ n } \zeta_k^{ 2 ( n - 1 ) } \exp\Bigg( - \binom{n}{2} \| \bfx_0 \|_1 \Bigg).
\end{align*}
By Stirling's formula and $\log k \geq n^2 / C^2$,
\begin{align*}
&\Pi^{ K, n }( \bfx \in ( 0, \infty )^{ n - 1 } : x_1 \geq \tilde f_n, \|\bfx - \bfx_0|_n \|_2 \leq  B_n^{-1} \tilde f^{nw}_n \zeta_k^2 ) \\
&\geq 2 \pi ( \log k )^{ 3 ( n - 1 ) / 2 } \\
&\phantom{ \geq 2 \pi  } \times  \exp\Bigg( 2n \log( n / e ) + ( n - 1 ) ( \log B_n' + \log C ) - \binom{n}{2} \| \bfx_0 \|_1  - ( n - 1 ) \log k \Bigg)  \\
&\sim2 \pi \exp\Bigg( 2n \log( n / e ) + ( n - 1 ) ( \log B_n' + \log C ) - \binom{n}{2} \| \bfx_0 \|_1  - ( n - 1 ) \frac{ n^2 }{ C^2 } \\
&\phantom{\sim2 \pi \exp\Bigg( 2n \log( n / e ) + ( n - 1 ) ( \log B_n' + \log C ) - \binom{n}{2} \| \bfx_0 \|_{ \infty } } + 3 ( n - 1 ) \log( n / C ) \Bigg).
\end{align*}
By Lemma \ref{lem:kingman_shape} we have that $\| \bfx_0 \|_1$ is bounded, and that $\log B_n'$ satisfies
\begin{align*}
w n \log n \geq \log B_n' &\geq w n \log \tilde f_n - \log A - \log \| \partial \psi^{ ( n ) } \|_{ \infty } - \frac{ 3 }{ 2 } \log n - n \log B \\
&\gtrsim - 4 w n \log n - \log A - \frac{ 5 }{ 2 } \log n - n \log B
\end{align*}
because $\tilde{ f }_n \geq \min_{ i \in \{ 1, \ldots, n - 1 \} }\{ x_{ 0, i } \} / 2$ and $\| \partial \psi^{ ( n ) } \|_{ \infty } \leq n$.
Hence, for sufficiently large $n$,
\begin{equation*}
\Pi^{ K, n }( \bfx \in ( 0, \infty )^{ n - 1 } : x_1 \geq \tilde f_n, \|\bfx - \bfx_0|_n \|_2 \leq  B_n^{-1} \tilde f^{nw}_n \zeta_k^2 ) \geq \exp\Big( - ( 1 + o( n ) ) \frac{ n^3 }{ C^2 } \Big).
\end{equation*}
Condition \eqref{eq:prior1} follows because
\begin{equation*}
\frac{ n^3 }{ C^2 } = C ( \log k )^{ 3 / 2 } = k \zeta_k^2.
\end{equation*}
Hence Theorem \ref{thm:general_rooted} holds, which shows that the sequence $\Pi^{K, n}$ is consistent at the claimed rate.

\subsection{Proof of Corollary \ref{thm:unconstrained_main}}

Before proving Corollary \ref{thm:unconstrained_main}, we collect some useful properties of the infinite-leaf uniform tree $( T, \bfx ) \sim \Pi^U$ in a preparatory lemma.

\begin{lem}\label{lem:uniform_shape}
With $\Pi^U$-probability one, $\lambda \| \bfx |_n \|_1 / n \to 1$ as $n \to \infty$, and only finitely many of the events $\{ x_{ n - 1 } < n^{ -2 } \}_{ n \geq 2 }$ happen.
\end{lem}
\begin{proof}
The first claim is just the strong law of large numbers for independent, $\text{Exp}(\lambda)$-distributed random variables.

For the second claim, Lemma \ref{lem:branches_bounds} yields
\begin{equation*}
\sum_{ n = 2 }^{ \infty }\Pi^U( x_{ n - 1 } < n^{ -2 } ) \leq \sum_{ n = 2 }^{ \infty }  \frac{ \lambda }{ n^2 } \Big( 1 + \frac{ \lambda }{ 2 n^2 } \Big) = \sum_{ n = 1 }^{ \infty } \frac{ \lambda ( 2 n^2 + \lambda ) }{ 2 n^4 } < \infty.
\end{equation*}
The claim follows by the Borel--Cantelli lemma.
\end{proof}

We will follow the same argument as the proof of Corollary \ref{thm:main} by checking \eqref{eq:entropy1}, \eqref{eq:remaining1}, and  \eqref{eq:prior1}, and concluding the claimed convergence using Theorem \ref{thm:general_rooted}.
For $\beta > 1$, $\delta > 0$, and $C > 0$ as in Corollary \ref{thm:unconstrained_main}, define
\begin{align*}
\zeta_k^2 &= \frac{ 2 C ( \log k )^{ 3 / 2 } }{ k }, \\
\eps_k^2 &= \frac{ 2 C ( \log k )^{ 3 / 2 + \beta } }{ k },\\
f_k &= \frac{ 1 }{ 2 n \lambda } \exp( - ( C + 4 ) k \zeta_k^2 ),\\
g_k &= \frac{ ( C + 4 + \delta ) }{ \lambda } k \zeta_k^2.
\end{align*}
Then, using $\log k \geq n^2 / C^2$,
\begin{align}
&\frac{ g_k^{ 2 n - 3 } }{\eps_k^{ 2 ( 2 n - 3 ) } f_k^{ w n ( 2 n - 3 ) } } \exp( - k \eps_k^2 ) \notag \\
&= \frac{ ( C + 4 + \delta )^{ 2 n - 3 } ( 2 n )^{ w n ( 2 n - 3 ) } \lambda^{ ( 2 n - 3 )( w n - 1 ) } }{ ( \log k )^{ \beta ( 2 n - 3 ) } } \notag \\
&\phantom{=} \times  \exp\Bigg( [ 2 n - 3 + 2 C ( C + 4 ) w n ( 2 n - 3 ) ( \log k )^{ 1 / 2 }  - 2 C ( \log k )^{ 1 / 2 + \beta } ] \log k \Bigg) \notag \\
&\leq \exp\Bigg( (2 n - 3 ) [ \log( C + 4 + \delta ) + ( w n - 1 ) \log \lambda - \beta \log ( n^2 / C^2 )  ] \notag \\
&\phantom{\sim \exp \Bigg( } + \Bigg[ w \frac{ ( 2 n - 3 ) }{ n } \log( 2 n ) + \frac{ 2 n - 3 }{ C^2 } + \frac{ 2 ( C + 4 ) }{ C^2 } w (2 n - 3) n^2  - \frac{ 2 n^{ 1 + 2 \beta } }{ C^{ 2 + 2 \beta } } \Bigg] n^2 \Bigg). \label{cor2_lhs}
\end{align}
Recalling \eqref{uniform_topologies}, the definition of $B_n$ in \eqref{Bn_definition}, and the fact that $\| \partial \psi \|_{ \infty } \equiv 1$, the right-hand side of \eqref{eq:entropy1} is bounded below as
\begin{align}
\frac{ 1 }{ | \TTT_n^U | ( 4 B_n )^{ n - 1 } } &\sim 2^{ n - 3 }\sqrt{ \frac{ n - 3 }{ 2n - 5 } } \frac{ e^{ n - 8 } ( n - 3 )^{ n - 3 } }{ ( 2 n - 5 )^{ 2 n - 5 } } \frac{ 1 }{ A \| \partial \psi \|_{ \infty } ( n - 1 ) B^n } \notag \\
&\geq \frac{  ( 1 + o( n ) )  }{ \sqrt{2} A n } \exp( ( n - 3 )[ \log 2 + \log( n - 3 ) ] + n - 8 \notag\\
&\phantom{\geq \frac{  ( 1 + o( n ) )  }{ \sqrt{2} A n } \exp(} - ( 2 n - 5 ) \log ( 2 n - 5 ) - n \log B - \log n ). \label{cor2_rhs}
\end{align}
When $\beta > 1$, \eqref{cor2_lhs} decays faster than \eqref{cor2_rhs} as $n \to \infty$.
Hence \eqref{eq:entropy1} holds for any sufficiently large $n$.

In order to show \eqref{eq:remaining1}, we use Lemma \ref{lem:branches_bounds}.
When $f_k\leq 1$
\[
\Pi^{ U, n }( \cup_{ T \in \TTT_n^U } L^c_{ f_k }( T )  )\leq n \lambda f_k \Big( 1 + \frac{ n \lambda }{ 2 } f_k \Big) \leq\frac{5}{8} \exp( - ( C + 4 ) k \zeta_k^2 ).
\]
On the other hand,
\begin{align*}
\Pi^{ U, n }( \cup_{ T \in \TTT_n^U } U^c_g( T ) ) \leq ( 2 n - 3 ) \exp( - \lambda g_k ) &= ( 2 n - 3 ) \exp( - ( C + 4 + \delta ) k \zeta_k^2 ) \\
&\leq \frac{ 3 }{ 8 } \exp( - ( C + 4 ) k \zeta_k^2 )
\end{align*}
when $k$ is large enough that $- \delta C ( \log k )^{ 3 / 2 } + \log( 2 n - 3 ) \leq \log ( 3 / 8 )$, which exists because $k$ grows faster than linearly with $n$.
Therefore, \eqref{eq:remaining1} is satisfied.

The left-hand side of \eqref{eq:prior1} does not depend on $T_0$ in this case because topologies and branch lengths are independent under $\Pi^{U, n}$ and $\Pi^U$, and the marginal distribution over topologies is uniform.
Let $B_n' := \tilde{f}_n^{nw} / (\sqrt{n}B_n)$.
We have by definition that $\| \partial \psi^{ ( n ) } \|_{ \infty } \equiv 1$ and  $\tilde{ f }_n \leq 1 / (2 \eta)$.
Hence
\begin{align}
B_n' \zeta_k^2 = \frac{ \tilde{ f }_n^{ n w } }{ \sqrt{ n } B_n } \frac{ 2 C ( \log k )^{ 3 / 2 } }{ k } &\leq \frac{1 }{ ( 2 \eta )^{ w n } n^{ 3/2 } A B^n } \frac{ 2 n^3 }{ C^2 \exp( n^2 / C^2 ) } \notag\\
&\sim \frac{ n^{ 3/2 } }{ A C^2 ( 2 \eta )^{ n w } B^n \exp( n^2 / C^2 ) } \to 0. \label{Bn_decay_uniform}
\end{align}
Again by Lemma \ref{lem:uniform_shape}, $\min\{ x_1, \ldots, x_{ n - 1 } \} \to 0$ no faster than $n^{-2}$, which is slower than the decay in \eqref{Bn_decay_uniform}.
Hence we can take a large enough $n$ and $k$ so that $x_{ 0, i } - B_n' \zeta_k^2 > 0$ for each $i \in \{ 1, \ldots, 2 n - 3 \}$, and $| x_i - x_{ 0, i } | \leq B_n' \zeta_k^2 \Rightarrow x_i \geq \tilde f_n$ for all $i \in \{ 1, \ldots, n \}$, whereupon
\begin{align*}
&\Pi^{ U, n }( \bfx \in ( 0, \infty )^{ 2 n - 3 } : \bfx \in L_{ \tilde{ f }_n }( T_0 ), \| \bfx - \bfx_0|_n \|_2 \leq B_n^{-1}\tilde f_n^{nw} \zeta_k^2 ) \\
&=\lambda^{ 2 n - 3 } \int_{\bfx \in \bbR_+^{2n-3}: \,\bfx \in L_{ \tilde{ f }_n }( T_0 ), \, \| \bfx - \bfx_0|_n \|_2 \leq \sqrt{n} B_n' \zeta_k^2 }\exp\left( - \lambda \sum_{ i = 1 }^{ 2 n - 3 } x_i \right) \rmd \bfx \\
&\geq \lambda^{ 2 n - 3 } \int_{ \bfx \in \bbR_+^{ 2 n - 3 } : \, | x_i - x_{ 0 , i } | \leq B_n' \zeta_k^2 }\exp\left( -\lambda \sum_{ i = 1 }^{ 2 n - 3 } x_i \right) \rmd \bfx \\
&=\prod_{ i = 1 }^{ 2 n - 3 } \left( \exp(-\lambda ( x_{ 0 , i } - B_n' \zeta_k^2 ) ) - \exp( - \lambda ( x_{ 0 , i } + B_n' \zeta_k^2 ) ) \right).
\end{align*}
We use a Taylor expansion to show that for any $i \in \{ 1, \dots, 2 n - 3 \},$
\[
\exp( - \lambda ( x_{ 0 , i } - B_n' \zeta_k^2 ) ) - \exp( - \lambda ( x_{ 0 , i } + B_n' \zeta_k^2 ) ) \geq 2 \lambda B_n' \zeta_k^2 \exp( - \lambda x_{ 0, i } ).
\]
Therefore, using $\log k \geq n^2 / C^2$ and $\| \bfx_0|_n \|_1 \sim n / \lambda$ from Lemma \ref{lem:uniform_shape},
\begin{align*}
&\Pi^{U, n}( \bfx \in ( 0, \infty )^{ 2 n - 3 } : \bfx \in L_{ \tilde{ f }_n }( T_0 ), \| \bfx - \bfx_0|_n \|_2 \leq \sqrt{n} B_n' \zeta_k^2 ) \\
&\geq ( 2 \lambda B_n' )^{ 2 n - 3 } \zeta_k^{ 2 ( 2 n - 3 ) } \exp( - \lambda \| \bfx_0|_n \|_1 ) \\
&\geq \frac{ \tilde{f}_n^{w n (2 n - 3) } ( 2C )^{ 2 n - 3 } ( \log k )^{ 3 ( 2 n - 3 ) / 2 } }{ A^{ 2 n - 3 } n^{3 ( 2 n - 3 ) / 2} B^{ n ( 2 n - 3 ) } k^{ 2 n - 3 } } \exp( - \lambda \| \bfx_0|_n \|_1 )  \\
&\geq \exp\Bigg( ( 2 n - 3 ) \Bigg[ w n \log \tilde{f}_n + \log (2 C) + 3 \log( n / C ) - \log A - \frac{ 3 }{ 2 } \log n - n \log B - \frac{ n^2 }{ C^2 } \Bigg] \\
&\phantom{\geq \exp\Bigg( } - (1 + o(n)) n \Bigg).
\end{align*}
By definition $\tilde{f}_n \leq 1 / ( 2 \eta )$.
To obtain a lower bound, note $\tilde{f}_n \geq \min\{ x_{ 0, 1 }, \ldots, x_{0, 2 n - 3 } \} / 2$.
By Lemma \ref{lem:uniform_shape}, $\min\{ x_{ 0, 1 }, \ldots, x_{0, 2 n - 3 } \} / 2 < n^{-2 } / 2$ only finitely often as $n \to \infty$.
These particularly short branches will all be internal eventually $\Pi^U$-almost surely by the second Borel--Cantelli lemma, because any given branch gets chosen as the attachment point of a new leaf with probability $1 / (2 n - 3)$ when there are $n$ leaves, and successive attachment points are independent.
Hence $\tilde{f}_n > n^{-2} / 2$ for any sufficiently large $n$, and
\begin{align*}
&\Pi^{U, n}( \bfx \in ( 0, \infty )^{ 2 n - 3 } : \bfx \in L_{ \tilde{ f }_n }( T_0 ), \| \bfx - \bfx_0|_n \|_2 \leq \sqrt{n} B_n' \zeta_k^2 ) \\
&\geq \exp\Bigg( ( 2 n - 3 ) \Bigg[ \log (2 C) + 3 \log( n / C ) - 2 w n \log ( n / \sqrt{2} )  - \log A - \frac{ 3 }{ 2 } \log n - n \log B \Bigg] \\
&\phantom{\geq \exp\Bigg( } - (1 + o(n)) n - \frac{ 2 n - 3 }{ 2 n } \frac{ 2 n^3 }{ C^2 } \Bigg) \\
&= \exp\Bigg( -(1 + o(n) ) \frac{ 2 n^3 }{ C^2 } \Bigg).
\end{align*}
Condition \eqref{eq:prior1} follows because
\begin{equation*}
\frac{ 2 n^3 }{ C^2 } = 2 C ( \log k )^{ 3 / 2 } = k \zeta_k^2.
\end{equation*}
Hence Theorem \ref{thm:general_rooted} holds, which shows that the sequence $\Pi^{U, n}$ is consistent at the claimed rate.

\section{Proofs of technical lemmas}\label{section:technical_proofs}

\subsection{Proof of Lemma \ref{lem:distances}}

For $i \in \{1, \ldots, m( T ) \}$, we shorten $\partial_i \psi_T( \bfx )_{ ( u, v ) }$ for the partial derivative $[ \partial \psi_T( \bfx ) / ( \partial x_i ) ]_{ ( u, v ) }$.
Then
\begin{align*}
\Bigg| \frac{ \partial P_{ ( T, \bfx ) }( \bfa ) }{ \partial x_i } \Bigg| &\leq \sum_{ \bfa^* \in H^{ | I( T ) | } }  r( a^{ ( u_0 ) } ) \sum_{ ( u, v ) \in T } | ( \partial_i \psi_T( \bfx )_{ ( u, v ) } \bfQ \exp( \bfQ \psi_T( \bfx )_{ ( u, v ) } ) )_{ a^{ ( u ) } a^{ ( v ) } } |\\
&\phantom{\leq \sum_{ \bfa^* \in H^{ | I( T ) | } }  r( a^{ ( u_0 ) } ) \sum_{ ( u, v ) \in T }} \times  \prod_{ ( u', v' ) \neq ( u, v ) }  \exp( \bfQ \psi_T( \bfx )_{ ( u', v' ) } )_{ a^{ ( u' ) } a^{ ( v' ) } }  \\
&\leq \sum_{ \bfa^* \in H^{ | I( T ) | } } \sum_{ ( u, v ) \in T : \partial_i \psi_T( \bfx ) \neq 0 } \| \partial \psi_T( \cdot ) \|_{ \infty } | ( \bfQ \exp( \bfQ \psi_T( \bfx )_{ ( u, v ) } ))_{ a^{ ( u ) } a^{ ( v ) } } | \\
&\leq \gamma \| \partial \psi^{ ( n ) } \|_{ \infty }( 2 n - 2 ) | H |^{ n - 1 },
\end{align*}
where we have bounded probabilities by 1, and used the fact that a tree with $n$ leaves has at most $n - 1$ internal nodes and $2 n - 2$ branches.
On the event $L_f( T )$ we also have
\begin{equation*}
P_{ ( T, \bfx ) }( \bfa ) \geq \min_{ \alpha \in H }\Big\{ r( \alpha ) \min_{ t \geq 0 } \{ \exp( \bfQ t )_{ \alpha \alpha } \}^{ n - 2 } \min_{ \beta \in H, t \geq f }\{ \exp( \bfQ t )_{ \alpha \beta } \}^n \Big\} =: \lambda_f,
\end{equation*}
where the minima are guaranteed to be strictly positive because $\exp( \bfQ t )$ has a stationary distribution, and because all external edges have strictly positive minimum lengths.
Hence
\begin{equation*}
\Bigg| \frac{ \partial \log P_{ ( T, \bfx )}( \bfa ) }{ \partial x_i } \Bigg| \leq \frac{ \gamma \| \partial \psi^{ ( n ) } \|_{ \infty }( 2 n - 2 ) | H |^{ n - 1 } }{ \lambda_f },
\end{equation*}
and therefore, whenever $\bfx_0 \in L_f( T_0 )$ and $\bfx \in L_f( T_0 )$ for $f > 0$,
\begin{equation*}
\KL( T_0, \bfx_0; T_0, \bfx ) \leq \frac{ \gamma \| \partial \psi^{ ( n ) } \|_{ \infty }( 2 n - 2 ) | H |^{ n - 1 }  }{ \lambda_f } \| \bfx_0 - \bfx \|_2.
\end{equation*}
It remains to bound $\lambda_f$ from below.
By definition, $w$ mutations suffice for changing a type to any other along a lineage with positive probability.
Moreover, mutations are distributed according to a Poisson process with rate at least $\eta > 0$.
For a sufficiently small branch length $f > 0$, we thus have
\begin{align*}
\lambda_f &\geq \min_{ \alpha \in H }\Big\{ r( \alpha ) \min_{ t \geq 0 } \{ \exp( \bfQ t )_{ \alpha \alpha } \}^{ n - 2 } \Big\} \frac{ ( \eta f )^{ n w } \exp( - n f ) }{ w! } \\
&\geq \min_{ \alpha \in H }\Big\{ r( \alpha ) \min_{ t \geq 0 } \{ \exp( \bfQ t )_{ \alpha \alpha } \}^{ n - 2 } \Big\} \frac{ ( \eta f )^{ n w } ( 1 - \eta f ) }{ w! },
\end{align*}
where the right hand side on the first line is a lower bound for the probability that the mutation process produces the observed configuration of types at the leaves.
It is obtained as the probability that type $\alpha$ appears at each internal node (of which there are at most $n - 2$ and the root), and that a sufficient number of mutations happens on each external branch (all of which have length at least $f$) to yield the desired type at each leaf.
Hence
\begin{align*}
&\KL( T_0, \bfx_0; T_0, \bfx ) \\
&\leq \max_{ \alpha \in H } \Big\{ \frac{ 1 }{ r( \alpha ) \min_{ t \geq 0 } \{ \exp( \bfQ t )_{ \alpha \alpha } \}^{ n - 2 } } \Big\} \frac{  \gamma \| \partial \psi^{ ( n ) } \|_{ \infty }( 2 n - 2 ) | H |^{ n - 1 } w! }{ ( \eta f )^{ n w } ( 1 - \eta f ) } \| \bfx_0 - \bfx \|_2,
\end{align*}
whereupon using $f < 1 /  (2 \eta)$ ensures that the right hand side is positive and yields the required bound.
The calculation for $V_0$ is essentially identical, and omitted.

\subsection{Proof of Lemma \ref{lem:branches_bounds}}

By Taylor's theorem with Lagrange form remainder,
\begin{equation*}
\bbP( X_i < f ) = 1 - \exp( - \lambda_i f ) \leq \lambda_i f + \frac{ \lambda_i^2 }{ 2 } \exp( - \lambda_i \zeta ) f^2
\end{equation*}
for some $\zeta \in ( 0, f )$.
Maximizing the right hand side with $\zeta = 0$ gives the claimed bound.

On the other hand, by the union bound,
\begin{align*}
\bbP( \max\{ X_1, \ldots, X_n \} > g ) &\leq \sum_{ i = 1 }^n \bbP( X_i > g ) = \sum_{ i = 1 }^n \exp( - \lambda_i g ) \leq n \exp( - \min\{ \lambda_1, \ldots, \lambda_n \} g ).
\end{align*}

\bibliographystyle{plainnat}
\bibliography{bibliography}

\begin{thebibliography}{31}
\providecommand{\natexlab}[1]{#1}
\providecommand{\url}[1]{\texttt{#1}}
\expandafter\ifx\csname urlstyle\endcsname\relax
  \providecommand{\doi}[1]{doi: #1}\else
  \providecommand{\doi}{doi: \begingroup \urlstyle{rm}\Url}\fi

\bibitem[Aldous(1993)]{aldous:1993}
D.~Aldous.
\newblock The continuum random tree {III}.
\newblock \emph{Ann. Probab.}, 21:\penalty0 248--289, 1993.

\bibitem[Borges and Kosiol(2020)]{borges/kosiol:2020}
R.~Borges and C.~Kosiol.
\newblock Consistency and identifiability of the polymorphism-aware
  phylogenetic models.
\newblock \emph{J. Theor. Biol.}, 486:\penalty0 110074, 2020.

\bibitem[Bouckaert et~al.(2019)Bouckaert, Vaughan, Barido-Sottani, Duch\^ene,
  Fourment, Gavryushkina, Heled, Jones, K\"uhnert, De~Maio, Matschiner, Mendes,
  M\"uller, Ogilvie, de~Plessis, Popinga, Rambaut, Rasmussen, Siveroni,
  Suchard, Wu, Xie, Zhang, Stadler, and Dummond]{bouckaertetal:2019}
R.~Bouckaert, T.~G. Vaughan, J.~Barido-Sottani, S.~Duch\^ene, M.~Fourment,
  A.~Gavryushkina, J.~Heled, G.~Jones, D~K\"uhnert, N.~De~Maio, M.~Matschiner,
  F.~K. Mendes, N.~F. M\"uller, H.~A. Ogilvie, L.~de~Plessis, A.~Popinga,
  A.~Rambaut, D.~Rasmussen, I.~Siveroni, M.~A. Suchard, C-H Wu, D.~Xie,
  C.~Zhang, T.~Stadler, and A.~J. Dummond.
\newblock {BEAST} 2.5: an advanced software platform for {Bayesian}
  evolutionary analysis.
\newblock \emph{PLOS Comput. Biol.}, 15:\penalty0 e1006650, 2019.

\bibitem[Dinh et~al.(2018)Dinh, {Si Tung Ho}, Suchard, and {Matsen
  IV}]{dinhetal:2018}
V.~Dinh, Lam. {Si Tung Ho}, M.~A. Suchard, and F.~A. {Matsen IV}.
\newblock Consistency and convergence rate of phylogenetic inference via
  regularization.
\newblock \emph{Ann. Stat.}, 46:\penalty0 1481--1512, 2018.

\bibitem[Felsenstein(1978)]{felsenstein:1978}
J.~Felsenstein.
\newblock The number of evolutionary trees.
\newblock \emph{Syst. Biol.}, 27:\penalty0 27--33, 1978.

\bibitem[Felsenstein(1981)]{felsenstein:1981}
J.~Felsenstein.
\newblock Evolutionary trees from {DNA} sequences: a maximum likelihood
  approach.
\newblock \emph{J. Mol. Evol.}, 17:\penalty0 368--376, 1981.

\bibitem[Freedman and Diaconis(1983)]{freedman/diaconis:1983}
D.~Freedman and P.~Diaconis.
\newblock On inconsistent {Bayes} estimates in the discrete case.
\newblock \emph{Ann. Stat.}, 11:\penalty0 1109--1118, 1983.

\bibitem[Ghosal and van~der Vaart(2017)]{gvdv}
S.~Ghosal and A.~W. van~der Vaart.
\newblock \emph{Fundamentals of Nonparametric Bayesian Inference}.
\newblock Cambridge Series in Statistical and Probabilistic Mathematics.
  Cambridge University Press, 2017.
\newblock \doi{10.1017/9781139029834}.

\bibitem[Ghosal et~al.(2000)Ghosal, Ghosh, and van~der Vaart]{ghosaletal:2000}
S.~Ghosal, J.~K. Ghosh, and A.~W. van~der Vaart.
\newblock Convergence rates of posterior distributions.
\newblock \emph{Ann. Stat.}, 28:\penalty0 500--531, 2000.

\bibitem[H{\"o}hna et~al.(2016)H{\"o}hna, Landis, Heath, Boussau, Lartillot,
  Moore, Huelsenbeck, and Ronquist]{hohna2016revbayes}
S.~H{\"o}hna, M.J. Landis, T.A. Heath, B.~Boussau, N.~Lartillot, B.R. Moore,
  J.P. Huelsenbeck, and F.~Ronquist.
\newblock {RevBayes: Bayesian phylogenetic inference using graphical models and
  an interactive model-specification language}.
\newblock \emph{Syst. Biol.}, 65\penalty0 (4):\penalty0 726--736, 2016.

\bibitem[Kingman(1982)]{Kingman82a}
J.~F.~C. Kingman.
\newblock The coalescent.
\newblock \emph{Stochast. Process. Applic.}, 13\penalty0 (3):\penalty0
  235--248, 1982.

\bibitem[Markin and Eulenstein(2021)]{markiv/eulenstein:2021}
A.~Markin and O.~Eulenstein.
\newblock Quartet-based inference is statistically consistent under the unified
  duplication-loss-coalescence model.
\newblock \emph{Bioinformatics}, 28:\penalty0 4064--4074, 2021.

\bibitem[Norris(1997)]{norris:1997}
J.~Norris.
\newblock \emph{Markov Chains}.
\newblock Cambridge University Press, 1997.

\bibitem[Paradis and Schliep(2019)]{paradis2019ape}
E.~Paradis and K.~Schliep.
\newblock ape 5.0: an environment for modern phylogenetics and evolutionary
  analyses in {R}.
\newblock \emph{Bioinformatics}, 35:\penalty0 526--528, 2019.
\newblock \doi{10.1093/bioinformatics/bty633}.

\bibitem[{R~Core Team}(2024)]{r2024}
{R~Core Team}.
\newblock \emph{{R: A Language and Environment for Statistical Computing}}.
\newblock R Foundation for Statistical Computing, Vienna, Austria, 2024.
\newblock URL \url{https://www.R-project.org/}.

\bibitem[Roch and Sly(2017)]{roch/sly:2017}
S.~Roch and A.~Sly.
\newblock Phase transition in the sample complexity of likelihood-based
  phylogeny inference.
\newblock \emph{Probab. Theory Related Fields}, 169:\penalty0 3--62, 2017.

\bibitem[Ronquist et~al.(2012)Ronquist, Teslenko, {van der Mark}, Ayres,
  Darling, H\"ohna, Larget, Liu, Suchard, and Huelsenbeck]{ronquistetal:2012}
F.~Ronquist, M.~Teslenko, P.~{van der Mark}, D.~L. Ayres, A.~Darling,
  S.~H\"ohna, B.~Larget, L.~Liu, M.~A. Suchard, and J.~P. Huelsenbeck.
\newblock Efficient {Bayesian} phylogenetic inference and model selection
  across a large model space.
\newblock \emph{Syst. Biol.}, 61:\penalty0 539--542, 2012.

\bibitem[Sayyari and Mirarab(2016)]{sayyari/mirarab:2016}
E.~Sayyari and S.~Mirarab.
\newblock Anchoring quartet-based phylogenetic distances and applications to
  species tree reconstruction.
\newblock \emph{BMC Genom.}, 17:\penalty0 101--113, 2016.

\bibitem[Schliep(2011)]{schliep2011phangorn}
K.P. Schliep.
\newblock phangorn: phylogenetic analysis in r.
\newblock \emph{Bioinformatics}, 27\penalty0 (4):\penalty0 592--593, 2011.
\newblock \doi{10.1093/bioinformatics/btq706}.

\bibitem[Schwartz(1965)]{schwartz:1965}
L.~Schwartz.
\newblock On {Bayes} procedures.
\newblock \emph{Z. Warsch. verw. Gebiete}, 4:\penalty0 10--26, 1965.

\bibitem[Semple and Steel(2013)]{semple/steel:2013}
C.~Semple and M.~A. Steel.
\newblock \emph{Phylogenetics}.
\newblock Oxford University Press, 2013.

\bibitem[Stamatakis(2006)]{stamatakis:2006}
A.~Stamatakis.
\newblock {RAxML-VI-HPC}: maximum likelihood-based phylogenetic analyses with
  thousands of taxa and mixed models.
\newblock \emph{Bioinformatics}, 22:\penalty0 2688--2690, 2006.

\bibitem[Steel(2013)]{steel:2013}
M.~A. Steel.
\newblock Consistency of {Bayesian} inference of resolved phylogenetic trees.
\newblock \emph{J. Theor. Biol.}, 336:\penalty0 246--249, 2013.

\bibitem[Steel and Sz\'ekely(2007)]{steel/szekely:2007}
M.~A. Steel and L.~A. Sz\'ekely.
\newblock Teasing apart two trees.
\newblock \emph{Comb. Probab. Comput.}, 16:\penalty0 903--922, 2007.

\bibitem[Steel and Sz\'ekely(2009)]{steel/szekely:2009}
M.~A. Steel and L.~A. Sz\'ekely.
\newblock Inverting random functions {III}: discrete {MLE} revisited.
\newblock \emph{Ann. Comb.}, 13:\penalty0 365--382, 2009.

\bibitem[Suchard et~al.(2018)Suchard, Lemey, Baele, Ayres, Drummond, and
  Rambaut]{suchardetal:2018}
M.~A. Suchard, P.~Lemey, G.~Baele, D.~L. Ayres, A.~J. Drummond, and A.~Rambaut.
\newblock Bayesian phylogenetic and phylodynamic data integration using {BEAST
  1.10}.
\newblock \emph{Virus Evol.}, 4:\penalty0 vey016, 2018.

\bibitem[Warnow(2017)]{warnow:2017}
T.~Warnow.
\newblock \emph{Computational Phylogenetics: An Introduction to Designing
  Methods for Phylogeny Estimation}.
\newblock Cambridge University Press, 2017.

\bibitem[Whidden and {Matsen IV}(2015)]{whidden/matsen:2015}
C.~Whidden and F.~A. {Matsen IV}.
\newblock Quantifying {MCMC} exploration of phylogenetic tree space.
\newblock \emph{Syst. Biol.}, 64:\penalty0 472--491, 2015.

\bibitem[Wickham(2016)]{wickham16}
H.~Wickham.
\newblock \emph{{ggplot2: Elegant Graphics for Data Analysis}}.
\newblock Springer-Verlag New York, 2016.

\bibitem[Zhang et~al.(2018)Zhang, Rabiee, Sayyari, and Mirarab]{zhangetal:2018}
C.~Zhang, M.~Rabiee, E.~Sayyari, and S.~Mirarab.
\newblock {ASTRAL-III}: polynomial time species tree reconstruction from
  partially resolved gene trees.
\newblock \emph{BMC Bioinform.}, 19:\penalty0 153, 2018.

\bibitem[Zhang et~al.(2021)Zhang, Dinh, and {Matsen IV}]{zhangetal:2021}
C.~Zhang, V.~Dinh, and F.~A. {Matsen IV}.
\newblock Nonbifurcating phylogenetic tree inference via the adaptive {LASSO}.
\newblock \emph{J. Am. Stat. Assoc.}, 116:\penalty0 858--873, 2021.

\end{thebibliography}

\begin{appendices}

\section{Examples}

This appendix contains additional details and figures for the experiments in Section~\ref{section:examples} of the main text.
The code to reproduce out experiments is hosted at \url{github.com/lukejkelly/TreeConsistency}.

\subsection{Experimental pipeline}

Our experimental pipeline for each tree prior is as follows.
\begin{description}
	\item[Trees:] We started from an initial tree with \( n = 4 \) leaves sampled from the prior. We then constructed a sequence of trees \( (T_0, \bfx_0) |_n \) on \( n = 4, 5, \dotsc, 16 \) leaves, each marginally distributed according to the prior, using the algorithms in Sections~4.1 and~4.2 of the main text.
	These trees in NEXUS format are available on request, instructions to generate data on them are included with the code in our GitHub repository.
	\item[Raw data:] For our experiments we consider trees with \( n = 4, 7, 10, 13, 16 \) leaves. On each tree \( (T_0, \bfx_0) |_n \), we generated 100 replicate data sets at \( 10^5 \) sites from a continuous-time Markov model with binary alleles and symmetric mutation rates \( Q_{0, 1} = Q_{1, 0} = \mu \) for each \( \mu = 0.0125, 0.025, 0.05, 0.1, 0.2 \).
	\item[Process data:] We created data sets of length \( k = 1, 10, \dotsc, 10^5 \) using the first \( k \) sites from each raw data set; that is, for each number of leaves \( n \), mutation rate \( \mu \) and replicate index \( r = 1, 2, \dotsc, 100 \), we formed a new data set using the alleles at sites \( 1, 2, \dotsc, k \) in the raw data set.
	For each tree prior, there are \( 15000 \) data sets since we consider all combinations of the number of leaves \( n \) (five values), sequence length \( k \) (six values), mutation rate \( \mu \) (five values) and replicate index (100 values).
	\item[Analyses:] We used Markov chain Monte Carlo for inference on each data set.
	Having fixed the mutation rate at its true value, we constructed a Markov chain targeting the posterior distribution of the tree given the sequence data.
	We used a standard suite of proposal mechanisms in the Markov chains.
	Each chain was initialized randomly and ran for two million iterations.
	We discarded the first $ 10^6 $ iterations as burn-in, then from the remainder drew $ 10^3 $ samples at intervals of $ 10^3 $ iterations.
	Mixing and convergence was diagnosed from trace plots of the log-likelihood at each iteration.
	\item[Figures:] For each set of samples, we computed a Monte Carlo estimate of the posterior support for the true tree topology.
	We then plot these estimates to illustrate how the posterior behaves as the number of leaves \( n \), number of sites \( k \) and mutation rate \( \mu \) vary across replicate data sets.
\end{description}
README files accompanying the code describe how we implement our pipeline in practice.

\subsection{Posterior support for true tree topology}
\label{app:posterior-support}

\begin{figure}
	\centering
	\includegraphics[width=\textwidth, trim=0cm 0cm 2cm 0cm, clip]{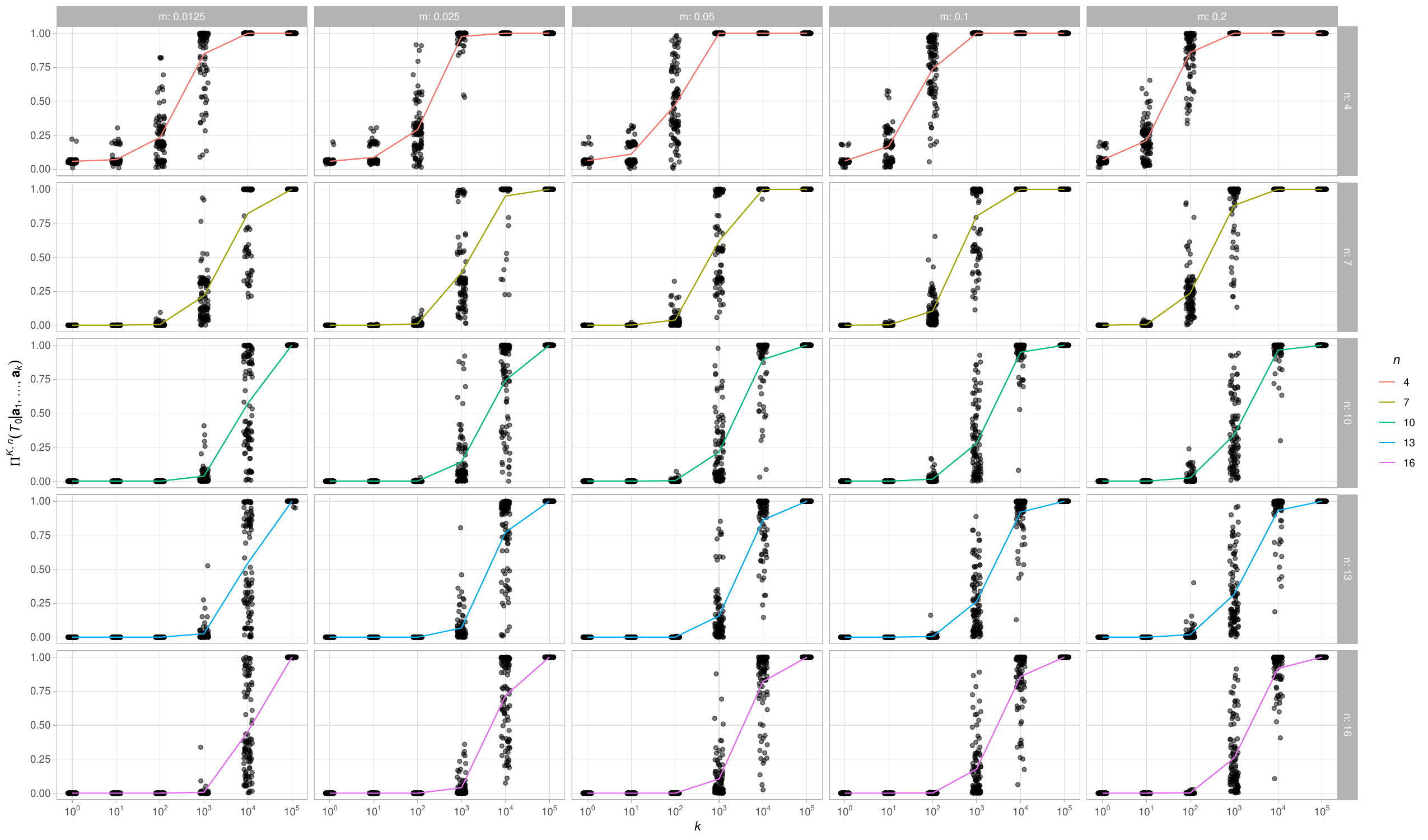}
	\caption{
		Posterior support for the true tree topology \( T_0 |_ n \) across 100 replicate data sets for each \( n \), \( \mu \) and \( k \) in experiments where the tree prior $ \Pi^{K, n} $ is Kingman's coalescent.
		Points depict the posterior support in a data set and curves display the average across replicate data sets.
		The curves are superimposed in \figref{fig:kingman-support} of the main text.
	}
	\label{fig:kingman-support-all}
\end{figure}

\begin{figure}
	\centering
	\includegraphics[width=\textwidth, trim=0cm 0cm 2cm 0cm, clip]{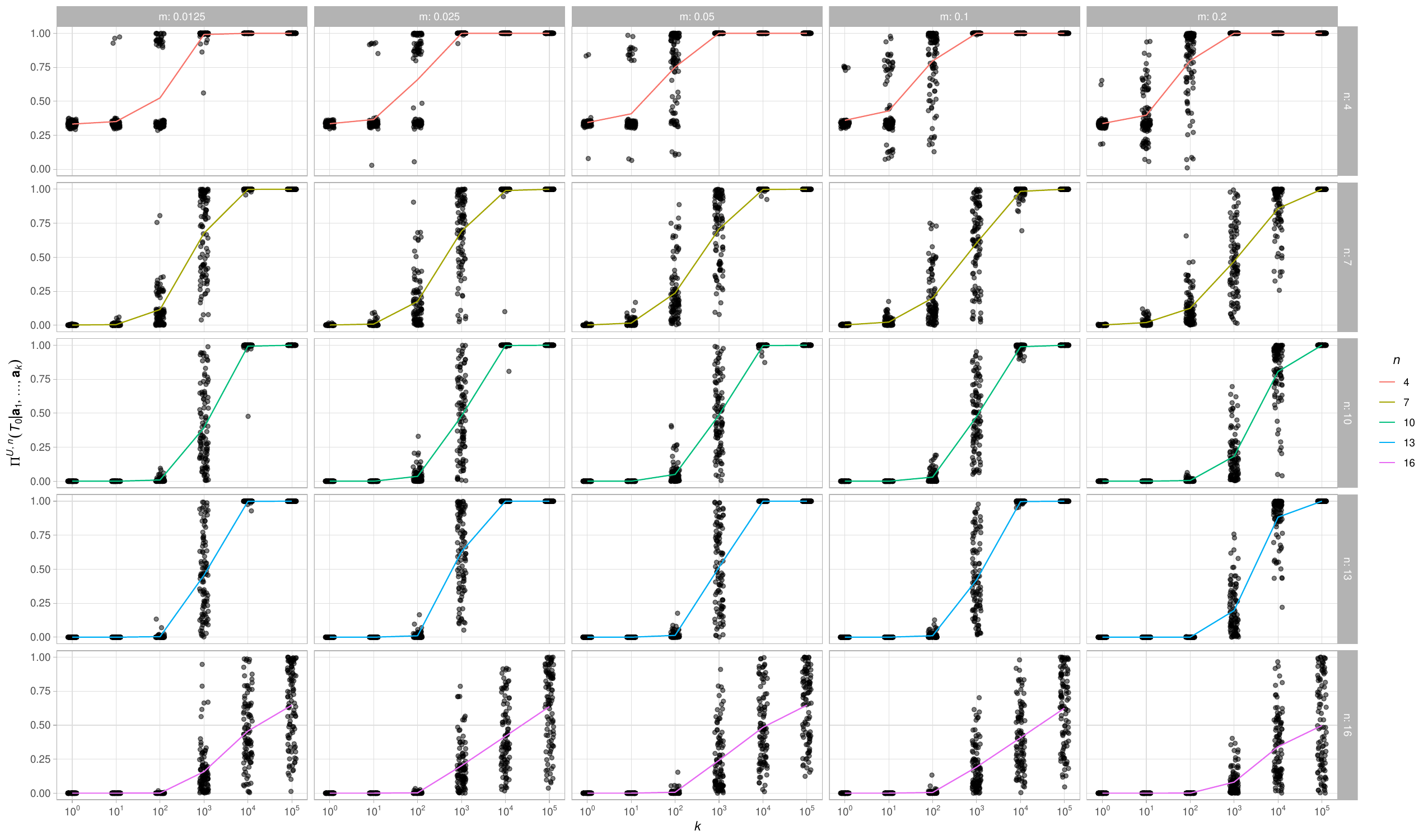}
	\caption{
		Posterior support for the true tree topology \( T_0 |_ n \) across 100 replicate data sets for each \( n \), \( \mu \) and \( k \) in experiments where the tree prior $ \Pi^{U, n} $ is uniform across topologies and exponential on branch lengths.
		Points depict the posterior support in a data set and curves display the average across replicate data sets.
		The curves are superimposed in \figref{fig:uniform-support} of the main text.
	}
	\label{fig:uniform-support-all}
\end{figure}

\figref{fig:kingman-support-all} displays the posterior support on the true rooted tree topology $ T_0 |_n $ under a Kingman coalescent prior $ \Pi^{K, n} $ in each replicate data set as $ n $, $ k $ and $ \mu $ vary.
In each case, the curves representing the average support approach 1 as \( k \) increases.
\figref{fig:uniform-support-all} displays the posterior support on the true unrooted tree topology $ T_0 |_n $ under the prior $ \Pi^{U, n} $ which is uniform across topologies and exponential on branch lengths.
For \( n = 16 \), we see comparatively slow growth in posterior support as \( k \) increases, particularly for \( \mu = 0.2 \).

We observe similar behaviour in both sets of experiments.
For small \( k \), the posterior support for the true topology is approximately equal to the prior support.
As \( k \) grows, we observe variation in posterior support across the replicate data sets and the average support starts to increase steadily.
For sufficiently large \( k \), the average posterior support approaches 1 and the variation across replicate data sets reduces.

\subsection{Markov chain behaviour}

\begin{figure}
	\centering
	\begin{subfigure}[t]{\textwidth}
		\centering
		\includegraphics[width=\textwidth, trim=0cm 245.7cm 0cm 0.75cm, clip]{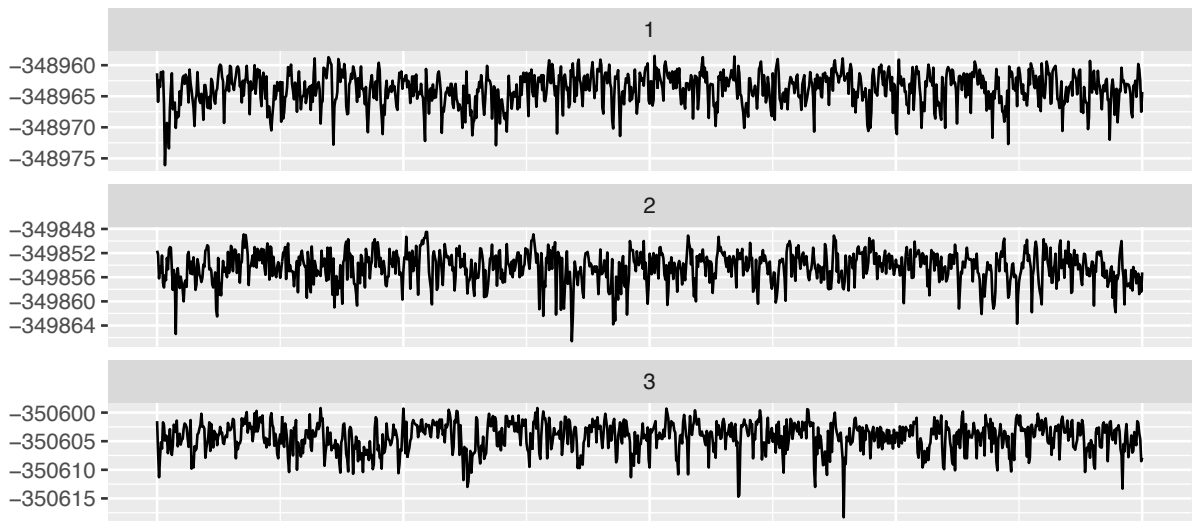}
		\caption{
			Kingman's coalescent prior.
		}
	\end{subfigure}

	\begin{subfigure}[t]{\textwidth}
		\centering
		\includegraphics[width=\textwidth, trim=0cm 245.7cm 0cm 0.75cm, clip]{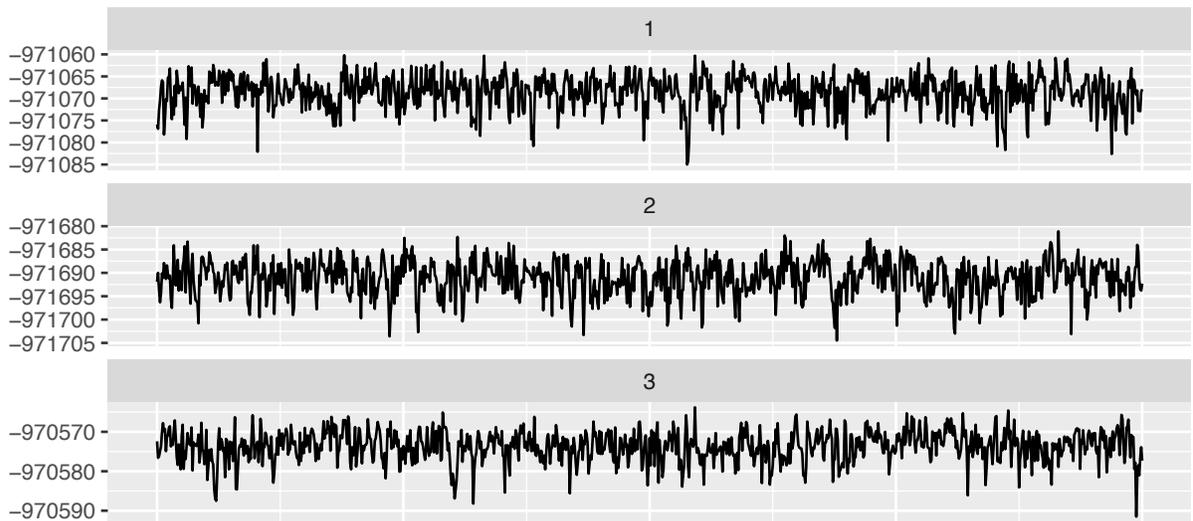}
		\caption{
			Tree prior is uniform across topologies and exponential on branch lengths.
		}
	\end{subfigure}
	\caption{
		Trace plots of the log-likelihood at each iteration in Markov chains used for inference in our experiments.
	}
	\label{fig:trace}
\end{figure}

Mixing and convergence of our Markov chains was monitored by visually inspecting trace plots of the log-likelihood at each iteration.
\figref{fig:trace} displays the log-likelihood after discarding burn-in for three data sets under each prior when \( n = 16 \), \( k = 10^5 \) and \( \mu = 0.2 \).
In each case, there is no visual indication that mixing or convergence failed to occur.

\end{appendices}

\end{document}